\documentclass[12pt]{amsart}

\usepackage{amssymb,latexsym}
\usepackage{graphicx,epsfig,color}

\def\cal{\mathcal}

\def\Bbb{\mathbb}

\def\ord{\text{\rm ord\,}}

\def\la{\langle}

\def\N{{\cal N}}
\def\T{{\cal T}}
\def\bC{{\Bbb C}}

\def\bN{{\Bbb N}}
\def\bR{{\Bbb R}}
\def\RR{{\Bbb R}}
\def\bZ{{\Bbb Z}}
\def\bQ{{\Bbb Q}}
\def\vp{{\varphi}}
\def\al{{\alpha}}
\def\be{{\beta}}
\def\ga{{\gamma}}
\def\la{{\lambda}}

\def\x{(x_1,x_2)}
\def\pa{{\partial}}
\def\all{{\rm all}}

\def\ka{{\kappa}}

\def\aol{{\left[ \begin{matrix} \al \\ l\end{matrix}\right ]}}

\def\bola{{\left[ \begin{matrix} \beta \\ \la\end{matrix}\right ]}}
\def\dotol{{\left[ \begin{matrix} \cdot \\ l\end{matrix}\right ]}}
\def\dotola{{\left[ \begin{matrix} \cdot \\ \la\end{matrix}\right ]}}

\def\bpm{\begin{pmatrix}}
\def\epm{\end{pmatrix}}

\def\noi{\noindent}
\def\bee{\begin{enumerate}}
\def\ee{\end{enumerate}}

\def\qed{\smallskip\hfill Q.E.D.\medskip}

\newtheorem{thm}{Theorem}[section]
\newtheorem{prop}[thm]{Proposition}
\newtheorem{cor}[thm]{Corollary}
\newtheorem{lemma}[thm]{Lemma}

\begin{document}

\title[Adapted coordinate systems]{On adapted coordinate systems}

\author[I. A. Ikromov]{Isroil A.Ikromov}
\address{Department of Mathematics, Samarkand State University,
University Boulevard 15, 703004, Samarkand, Uzbekistan}
\email{{\tt ikromov1@rambler.ru}}
\author[D. M\"uller]{Detlef M\"uller}
\address{Mathematisches Seminar, C.A.-Universit\"at Kiel,
Ludewig-Meyn-Stra\ss{}e 4, D-24098 Kiel, Germany} \email{{\tt
mueller@math.uni-kiel.de}}

\thanks{We acknowledge the support for this work be the Deutsche Forschungsgemeinschaft.}
\begin{abstract}
The notion of an adapted coordinate system, introduced by   V~.I. ~Arnol'd, plays an important role in the study of asymptotic expansions of oscillatory integrals.  In two dimensions, A.~N.~Varchenko gave sufficient conditions for the adaptness of a given coordinate system and proved the  existence of an adapted coordinate system for a class of analytic functions without  multiple components. Varchenko's  proof is based on Hironaka's theorem on the resolution of singularities.

In this article, we present a new, elementary and concrete approach to these results, which is based on the Puiseux  series expansion of roots of the given function.  Our method applies to arbitrary real analytic functions, and even extends to arbitrary smooth functions of finite type. Moreover, by avoiding  Hironaka's theorem, we can give necessary and sufficient conditions for the adaptedness of a given coordinate system in the smooth, finite type setting.

\end{abstract}
\maketitle

%\vfill\newpage

\tableofcontents

\thispagestyle{empty}

%\begin{document}
%\title{Estimates for damping Oscillatory integrals}
%\author{}
%\maketitle
\setcounter{equation}{0}
\section{Introduction}\label{intro}

It is an obvious  fact that the  asymptotic behavior of an oscillatory integral of the form
$$I(\la):=\int_{\bR^n} e^{i\la f(x)} a(x)\, dx
$$ 
does not change under  a smooth change of variables $x=\vp(y).$

 This observation is employed already in the proof of  van der Corput's  lemma (see, e.g., \cite{St}), according to which the   asymptotic behavior of a one-dimensional oscillatory
integral is determined by the maximal order of vanishing at the  critical points of the phase function $f.$ 

In higher dimensions, the problem  of determining the exact asymptotic behavior of an 
oscillatory integral is substantially more difficult. V.I. Arnol'd
 conjectured in \cite{A} that the asymptotic behavior of the  oscillatory
integral $I(\la)$  is determined by the Newton polyhedron associated to the phase function $f$ in a 
so-called  ``adapted''  coordinate system. For 
some special  cases this conjecture was then indeed verified by means of  Arnol'd's classification of singularities (see \cite{AGV}). Later, however,  A.N. Varchenko \cite{Va} disproved  Arnol'd's
conjecture in dimensions  three and higher.

Moreover, in the same paper he was able to verify Arnol'd's
conjecture for  two-dimensional oscillatory integrals. In particular, he proved in dimension two  the  existence of adapted coordinate systems, and showed   that the leading term of the 
asymptotic expansion of  $I(\la)$ can be  constructed from the Newton polyhedron associated to the phase function $f$ in such an adapted coordinate system.

The purpose of this article is to present a new, elementary and concrete approach to the latter  results in two dimensions, which is based on the Puiseux  series expansion of roots of the given function $f.$   Our method applies to arbitrary real analytic functions, and even extends to arbitrary smooth functions of finite type. Moreover, by avoiding  Hironaka's theorem, we can give necessary and sufficient conditions for the adaptedness of a given coordinate system in the smooth, finite type setting.

\bigskip

\setcounter{equation}{0}
\section{Preliminaries}\label{adaptdco}

Let $f$ be a smooth real-valued  function defined on a neighborhood of the origin in $\bR^2$ with $f(0,0)=0,\, \nabla f(0,0)=0,$ and consider the associated Taylor series 
$$f(x_1,x_2)\sim\sum_{j,k=0}^\infty c_{jk} x_1^j x_2^k$$
of $f$ centered at  the origin. 
The set
$$\T(f):=\{(j,k)\in\bN^2: c_{jk}=\frac 1{j!k!}\partial_{ x_1}^j\partial_{ x_2}^k f(0,0)\ne 0\}
$$ 
will be called the {\it Taylor support } of $f$ at $(0,0).$  We shall always assume that 
$$\T(f)\ne \emptyset,$$
i.e., that the function $f$ is of finite type at the origin. If $f$ is real analytic, so that the Taylor series converges to $f$ near the origin, this just means that $f\ne 0.$ 
The {\it Newton polyhedron} $\N(f)$ of $f$ at the origin is defined to be the convex hull of the union of all the quadrants $(j,k)+\bR^2_+$ in $\bR^2,$ with $(j,k)\in\T(f).$  The associated {\it Newton diagram}  $\N_d(f)$ in the sense of Varchenko \cite{Va}  is the union of all compact faces  of the Newton polyhedron; here, by a {\it face,} we shall  mean an edge or a vertex.

We shall use coordinates $(t_1,t_2)$ for points in the plane containing the Newton polyhedron, in order to distinguish this plane from the $(x_1,x_2)$ - plane. 

The {\it distance} $d=d(f)$  between the Newton
polyhedron and the origin in the sense of Varchenko is given by the coordinate $d$ of the point $(d,d)$ at which the bisectrix   $t_1=t_2$ intersects the boundary of the Newton polyhedron.

The {\it principal face} $\pi(f)$  of the Newton polyhedron of $f$  is the face of minimal dimension  containing the point $(d,d)$. Deviating from the notation in \cite{Va}, we shall call the series
$$f_p(x_1,x_2):=\sum_{(j,k)\in \pi(f)}c_{jk} x_1^j x_2^k
$$
the {\it principal part} of $f.$ In case that $\pi(f)$ is compact,  $f_p$ is a mixed homogeneous polynomial; otherwise, we shall  consider $f_p$ as a formal power series. 

Note that the distance between the Newton polyhedron and
the origin depends on the chosen local coordinate system in which $f$ is expressed.  By a  {\it local analytic (respectively smooth) coordinate system at the origin} we shall mean an analytic (respectively smooth)   coordinate system defined near the origin which preserves $0.$ If we work in the category of smooth functions $f,$  we shall always consider smooth coordinate systems, and if $f$ is analytic, then one usually restricts oneself to analytic coordinate systems (even though this will not really be necessary for the questions we are going to study, as we will see). The {\it height } of the analytic (respectively smooth) function $f$ is defined by 
$$h(f):=\sup\{d_x\},$$
 where the
supremum  is taken over all local analytic (respectively smooth) coordinate systems $x$ at the origin, and where $d_x$
is the distance between the Newton polyhedron and the origin in the
coordinates  $x$.

A given   coordinate system $x$ is said to be
 {\it adapted} to $f$ if $h(f)=d_x.$
 
A.N. Varchenko \cite{Va} proved that if $f$ is a real analytic function (without multiple components) near the origin in $\bR^2,$ then there exists a  local analytic coordinate system which is adapted to $f.$
The proof of this result in \cite{Va}  is based on Hironaka's deep theorem (see \cite{H}) on the resolution of singularities.

We shall here give an elementary proof of  Varchenko's 
theorem which is based on the Puiseux series expansion of roots. Moreover, 
our method extends to prove an analog of A.N. Varchenko
theorem for smooth functions.

It may be interesting at this point to remark that, of course, the notions introduced above extend to smooth functions in more than two real variables. However, as shown by Varchenko \cite{Va}, in dimensions higher than two adapted coordinate systems may not exist, even in the analytic setting.

\medskip
\subsection{The principal part of $f$ associated to a supporting line of the Newton polyhedron as a mixed homogeneous polynomial }

Let $\ka=(\ka_1,\ka_2)$ with $\ka_1,\ka_2>0$ be a given weight, with associated one-parameter family of dilations $\delta_r(x_1,x_2):= (r^{\ka_1}x_1,r^{\ka_2}x_2),\ r>0.$  A function  $f$ on $\bR^2$ is  said to be   {\it $\ka$-homogeneous of degree $a,$} if $f(\delta_rx)=r^a f(x)$ for every $r>0, x\in\bR^2.$ Such functions will also be called {\it mixed homogeneous.} 
The exponent $a$ will be denoted as the  {\it $\ka$-degree} of $f.$ For instance, the monomial $x_1^j x_2^k$ has $\ka$-degree $\ka_1 j+\ka_2k.$

If $f$ is an arbitrary smooth function near the origin,  consider its Taylor series $\sum_{j,k=0}^\infty c_{jk} x_1^j x_2^k$  around the origin.  We choose $a$ so that the line  $L_\ka:=\{(t_1,t_2)\in\bR^2:\ka_1 t_1+\ka_2t_2=a\}$ is the supporting  line to the Newton polyhedron $\N(f) $ of $f.$ Then the non-trivial polynomial 
$$f_\ka(x_1,x_2):=\sum_{(j,k)\in L_\ka} c_{jk} x_1^j x_2^k
$$
is $\ka$-homogeneous of degree $a;$ it will be called the {\it $\ka$-principal part} of $f.$  By definition, we then have 
\begin{equation}\label{princ1}
f(x_1,x_2)=f_\ka(x_1,x_2) +\ \mbox{terms of higher $\ka$-degree.}
\end{equation}
More precisely, we mean by this that  every point $(j,k)$ in the Taylor support of the   remainder term $f-f_\ka$ lies  on a line $\ka_1t_1+\ka_2t_2=d$ with $d>a$   parallel to, but above the line $L_\ka,$  i.e., 
we have $\ka_1j+\ka_2k>a.$ 
Moreover, clearly
$$
\N_d(f_\ka)\subset \N_d(f).
$$

The following lemma gives an equivalent description  of the notion ''terms of higher $\ka$-degree'', which is quite useful in applications, since it may be used to essentially reduce many considerations to the case of polynomial  functions. It will be mostly applied without further mentioning. 

\begin{lemma} \label{s2.1n} 
Assume that $f$ is a smooth function defined near the origin, and let $c\ge 0.$  Assume that $\ka_1\le \ka_2,$ and choose $m\ge 1$ in $\bN$ such that $\ka_1m >c.$ Then $f$ consists  of terms of $\ka$-degree greater or equal to $c$ in the   above sense (i.e., $\ka_1j+\ka_2k\ge c$ for every $(j,k)\in \T(f)$)  if and only if there exists a polynomial function $F$ with  
$ \quad \T(F)\subset\{(j,k)\in\bN:\ka_1 j+\ka_2 k\ge c\}$ and smooth functions 
$a_{jk},$  for $j+k=m,$ such that 
\begin{equation}\label{princ2}
f(x_1,x_2)=F(x_1,x_2) +\sum_{j+k=m}   x_1^jx_2^k a_{jk}\x .
\end{equation}
Notice that $\ka_1 j+\ka_2 k> c$ whenever $j+k=m.$ 
\end{lemma} 

\noi {\bf Proof.} Assume that $\ka_1j+\ka_2k\ge c$ for every $(j,k)\in \T(f).$ If  we then  choose for $F$ the Taylor polynomial of degree $m-1$ of $f$ at the origin,  the  representation \eqref{princ2} follows from Taylor's formula. 

Conversely, it is obvious that  the  representation \eqref{princ2} implies that $\ka_1j+\ka_2k\ge c$ for every $(j,k)\in \T(f).$

\qed

Sometimes, it will be convenient to extend these definitions to the case where $\ka_1=0$ or $\ka_2=0.$ 
In that case, the $\ka$- principal part $\sum_{(j,k)\in L_\ka} c_{jk} x_1^j x_2^k$ of $f$ will just be considered as a formal power series, unless $f$ is real analytic, when it is real analytic too. 

\medskip

Let $P\in \bR[x_1,x_2]$ be a mixed  homogeneous
polynomial,  and
assume that $\nabla P(0,0)=0$. Following  \cite{IKM}, we denote by 
$$m(P):=\ord_{S^1} P$$ 
the maximal order of vanishing of $P$ along the unit circle $S^1$ centered at the
origin. 

% Moreover, we define the {\it height of the homogeneous polynomial } $P$ by
%$$
%\tilde h(f):=max \{m(f), 1/|\ka|\},\quad \mbox{where}\quad
%|\ka|=\ka_1+\ka_2.
%$$

\medskip
The following Proposition will be a useful tool.

\medskip
If $m_1,\dots,m_n$ are positive integers, then we denote by
$(m_1,\dots,m_n)$ their greatest common divisor. 

\begin{prop}\label{s2.1} Let $P$ be a $(\ka_1,\ka_2)$-homogeneous
polynomial of degree one, and assume that $P$ is not of the form
$P(x_1,x_2)=cx_1^{\nu_1} x_2^{\nu_2}.$ Then $\ka_1$ and $\ka_2$ are
uniquely determined by $P,$ and  $\ka_1,\ka_2\in \bQ.$  

Let us assume that $\ka_1\le \ka_2,$ and write 
$$
\ka_1=\frac{q}{m},\, \ka_2=\frac{p}{m},\quad (p,q,m)=1,
$$
so that in particular $p\ge q.$
Then $(p,q)=1,$ and there exist non-negative integers $\al_1,\,\al_2$ and
a
$(1,1)$-homogeneous polynomial $Q$ such that the polynomial $P$ can be written as
\begin{equation}\label{2.3}
P(x_1,x_2)=x_1^{\al_1} x_2^{\al_2}Q(x_1^p,x_2^q).
\end{equation}
More precisely, $P$ can be written in the form
\begin{equation}\label{2.4}
P(x_1,x_2)=c x_1^{\nu_1}x_2^{\nu_2}\prod_{l=1}^M(x_2^q-\la_l
x_1^p)^{n_l},
\end{equation}
with $M\ge 1,$ distinct $\la_l\in \bC\setminus\{0\}$ and multiplicities $n_l\in \bN\setminus\{0\},$ with  $\nu_1,\nu_2\in\bN$ (possibly different from $\al_1,\,\al_2$  in
\eqref{2.3}). 

Let us   put  $n:=\sum_{l=1}^M n_l.$ The distance $d(P)$ of $P$ can then be read off from \eqref{2.4} as follows:

If the principal face of $\N(P)$ is compact, then it lies on the line $\ka_1 t_1+\ka_2 t_2=1,$ and the distance is given by 
\begin{equation}\label{2.8}
d(P)=\frac 1{\ka_1+\ka_2}=\frac{\nu_1 q+\nu_2 p+pqn}{q+p}.
\end{equation}
 Otherwise, we have $d(P)=\max\{\nu_1,\nu_2\}.$ In particular, in any case we have $d(P)=\max\{\nu_1,\nu_2,\frac 1{\ka_1+\ka_2}\}.$
\medskip

\end{prop}

\noindent{\bf Proof.} The proof of Proposition \ref{s2.1} is based on
elementary number  theoretic arguments. Denote by $A$ the set of all
solutions
$(\al,\beta)\in\bN^2$ of the linear equation
\begin{equation}\label{2.5}
\al \ka_1+\be\ka_2=1 .
\end{equation}

Then 
\begin{equation}\label{2.6}
P(x_1,x_2)=\sum_{(\al,\beta)\in A}c_{\al,\beta}\,  x_1^\al
x_2^\beta,
\end{equation}
for suitable coefficients $c_{\al,\beta}.$

 If $P$ is not of the form
$P(x_1,x_2)=cx_1^{\nu_1} x_2^{\nu_2},$ then the equation \eqref{2.5}  has
at least two different solutions $(\al_1,\beta_1), (\al_2,\beta_2)
\in\bN^2$ for which the coefficients
$c_{\al,\beta}$ in \eqref{2.6} are non-vanishing. The corresponding
equations $\al_1 \ka_1+\be_1\ka_2=1$ and  $\al_2 \ka_1+\be_2\ka_2=1$
determine the numbers $\ka_1,\ka_2$ uniquely and show that they are
rational.

\medskip

Assume next that $\ka_1\le \ka_2,$ and write 
$$
\ka_1=\frac{q}{m},\, \ka_2=\frac{p}{m},\quad (p,q,m)=1.
$$
Choose then  $(\al_0,\be_0)\in A$ with $\al_0$ maximal.
Since $p\bZ+q\bZ +m\bZ=\bZ,$  the identity
\eqref{2.5}, which is equivalent to 
$$\al q+\be p=m,
$$  implies that $p\bZ+q\bZ=\bZ$ (since $A\ne\emptyset$), so
that 
\begin{equation}\label{2.7}
(p,q)=1. 
\end{equation}

Notice that $(\al,\be)\in A$ is equivalent to
$(\al-\al_0)q+(\be-\be_0)p=0.$ Then \eqref{2.7} implies
that
$\al-\al_0=-sp, \, \be-\be_0=sq$ for some  $s\in \bZ$, or, equivalently,
$$
\al=\al_0-sp, \, \be=\be_0+sq. 
$$
Since $\al_0$ is  maximal  for $A$, we have $s\ge0$.
Choose $s_1\in\bZ $ maximal with  $\al_1:=\al_0-s_1p\ge 0.$ Then  for any
$(\al,\be)\in A$ we have the relation
$$
\al=\al_1+(s_1-s)p,\quad \be=\be_0+sq.
$$
Notice  that $s,\, s_1-s\in\bN.$ So, every monomial $x_1^\al
x_2^\be$ with $(\al,\be)\in A$ can be written as
$$
x_1^\al x_2^\be=x_1^{\al_1}x_2^{\be_0}(x_1^p)^{s_1-s}(x_2^q)^{s}.
$$
This in combination with \eqref{2.6} yields \eqref{2.3}.
\medskip

In order to prove \eqref{2.4},  write
$$Q(y_1,y_2)= cy_2^n +c_1 y_2^{n-1}y_1+\dots +c_n y_1^n,
$$
where $n$ is the degree of $Q.$ We may then assume that $c\ne 0,$ for
otherwise we can pull out some power of $y_1=x_1^p$ from $Q$ in
\eqref{2.3}. Assuming without loss of generality that $c=1,$ we can write 
$$Q(y_1,y_2)=y_1^n\, Q(1,\frac{y_2}{y_1})=y_1^n\prod_{j=1}^n
(\frac{y_2}{y_1}-\la_j)=\prod_{j=1}^n (y_2-\la_j y_1),
$$
where $\la_1,\dots,\la_n\in\bC$  are the roots of the polynomial $Q(1,y_2),$ listed with their multiplicities. This yields \eqref{2.4}. 
\medskip

To compute the distance $d(P),$ observe that the vertices of the Newton polyhedron of $P$ are given by
$ X_1:=(\nu_1,\nu_2+nq)$ and $X_2:= (\nu_1+np,\nu_2),$  which, by our assumptions,  are different points on the line $\ka_1t_1+\ka_2t_2=1.$ One then easily computes that 
$$\ka_1=\frac {q}{\nu_1q+\nu_2p+pqn},\ \ka_2=\frac {p}{\nu_1q+\nu_2p+pqn}.
$$
Thus, if the principal face of $\N(P)$ is compact, then it is the interval connecting these two points, hence it  lies on the line $\ka_1t_1+\ka_2t_2=1,$ and we immediately obtain \eqref{2.8}.  Otherwise, if  $\nu_1\le \nu_2, $ then the principal face is the horizontal half-line with left endpoint $X_1,$ so that $d(P)=\nu_2,$ and similarly $d(P)=\nu_1,$ if $\nu_1\ge \nu_2.$

From the geometry of $\N(P),$ it is then clear that we always have the identity $d(P)=\max\{\nu_1,\nu_2,\frac 1{\ka_1+\ka_2}\}.$

\qed

The proposition shows that every zero (or ``root'')  $\x$ of $P$ which does not lie on a coordinate axis is of the form $x_2=\la_l^{1/q} x_1^{p/q}.$ The quantity 
$$
d_h(P)=\frac 1{\ka_1+\ka_2}
$$
will be called the {\it homogeneous distance} of the mixed homogeneous polynomial $P.$ Recall that $(d_h(P),d_h(P))$ is just the point of intersection of the bisectrix with the line $\ka_1t_1+\ka_2t_2=1$ on which the Newton diagram $\N_d(P)$ lies. Moreover,
$$
d_h(P)\le d(P),
$$
and  the proof of Proposition \ref{s2.1} shows that 
$$ d_h(P)=\frac{\nu_1 q+\nu_2 p+pqn}{q+p}.$$
Notice also that 
\begin{equation}\label{2.9}
m(P)=\max\{\nu_1,\nu_2,\max_{l\in R} n_l\},
\end{equation}
where the index set  $R:=\{l=1,\dots,M:\la_l\in\RR\}$ corresponds to the set of real roots of $P$ which do not lie on a coordinate axis.

\begin{cor}\label{s2.3} 
 Let $P$ be a $(\ka_1,\ka_2)$-homogeneous
polynomial of degree one as in Proposition \ref{s2.1}, and consider the representation \eqref{2.4} of $P.$ We put again $n:=\sum_{l=1}^M n_l.$
\bee
\item[(a)]  If  $\ka_2/\ka_1\notin \bN,$ i.e., if $q\ge 2,$  then $n< d_h(P).$ In particular, every real root $x_2=\la_l^{1/q} x_1^{p/q}$ of $P$ has multiplicity $n_l< d_h(P).$

\item[(b)]  If $\ka_2/\ka_1\in \bN,$ i.e., if $q=1,$  then there exists at most one real root of $P$ on the unit circle $S^1$ of multiplicity greater than $d_h(P).$ More precisely, if  we put $n_0:=\nu_1, n_{M+1}:=\nu_2$ and choose $l_0\in\{0,\dots,M+1\}$ so that  $n_{l_0}=\max\limits_{l=0,\dots,M+1} n_l,$  then $n_l\le d_h(P)$ for every $l\ne l_0.$

\ee
\end{cor}

\noindent{\bf Proof.}  If  $q\ge 2,$ then $p>q\ge 2,$ so that $\frac 1p+\frac 1q<1,$ hence
$$ d_h(P)\ge \frac{pqn}{q+p}>n.
$$
This proves (a).

Assume next that $q=1,$ and suppose that $n_{l_1}>d_h(P)$ and  $n_{l_2}>d_h(P),$ where $l_1< l_2.$ If $l_1>0,$ then we arrive at the contradiction
$$d_h(P)\ge \frac {p(n_{l_1}+n_{l_2})}{p+1}>\frac{2pd_h(P)}{p+1}\ge d_h(P).
$$
Similarly, if $l_1=0,$ then we obtain the contradiction
$$d_h(P)\ge \frac {n_0+pn_{l_2}}{p+1}>\frac{d_h(P)+pd_h(P)}{p+1}=d_h(P).
$$

\qed

The corollary shows in particular that the multiplicity  of every real root of $P$ not lying  on a coordinate  axis is bounded 
by the distance $d(P),$ unless $q=1,$ in which case there can at most be one real root $x_2=\la_{l_0}x_1^p$ with multiplicity exceeding  $d(P).$ If such a root exists, we shall call it the {\it principal root of $P.$} 
\bigskip

\setcounter{equation}{0}
\section{Conditions for adaptedness of a given coordinate system}\label{conditions}

 In \cite{Va}, Varchenko has provided various sufficient conditions for the adaptedness of a given coordinate system (at least under certain non-degeneracy conditions), which prove to be useful. The goal of this section is to provide a new, elementary approach  to these results in the general case. Our approach has the additional advantage of extending to the category of smooth functions.
 
 \subsection{On the effect of a change of coordinates on the Newton diagram}\label{effect}
 We have to understand what effect a change of coordinates has on the Newton diagram. To this end, we shall make use of the following auxiliary result.
 % (in our corresponding  figures, $N_p(P)$ will denote the Newton polyhedron $\N(P)$).

\begin{figure}
\centering
\includegraphics[scale=0.7]{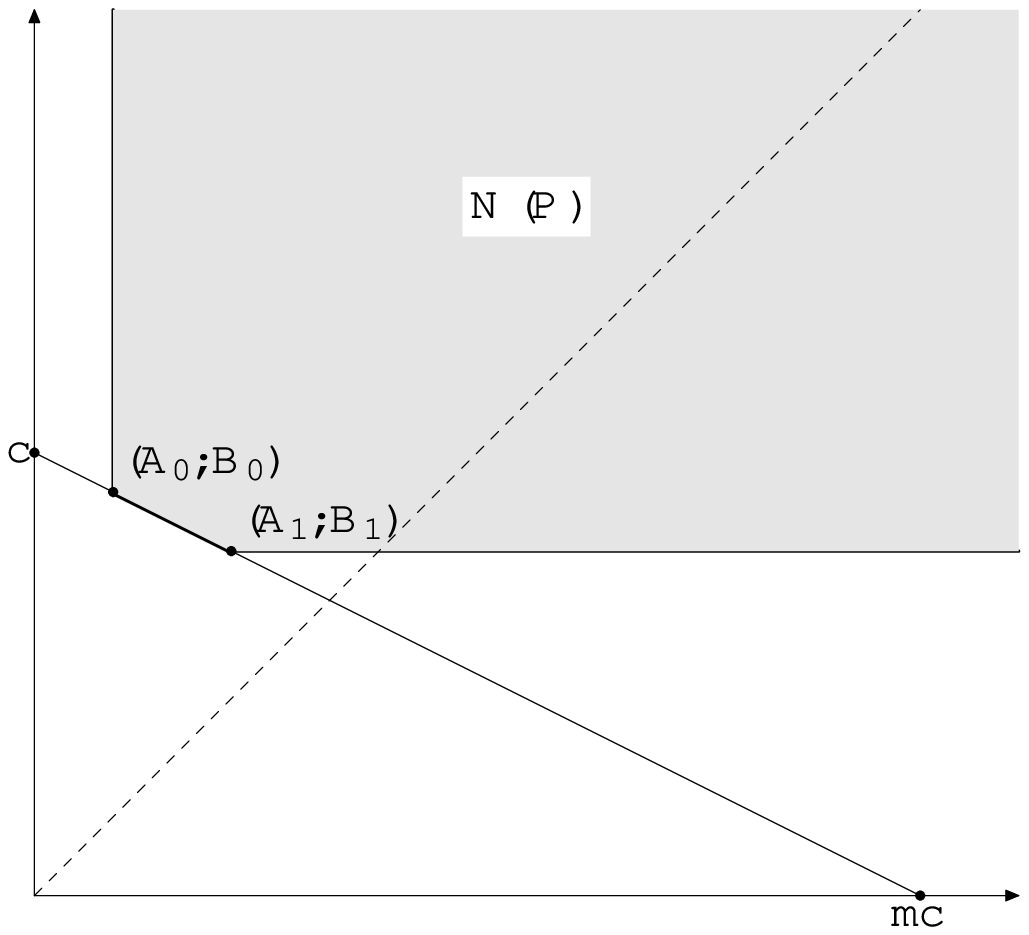}
\caption{\label{fig1}}
\end{figure}

\begin{figure}
\centering
\includegraphics[scale=0.7]{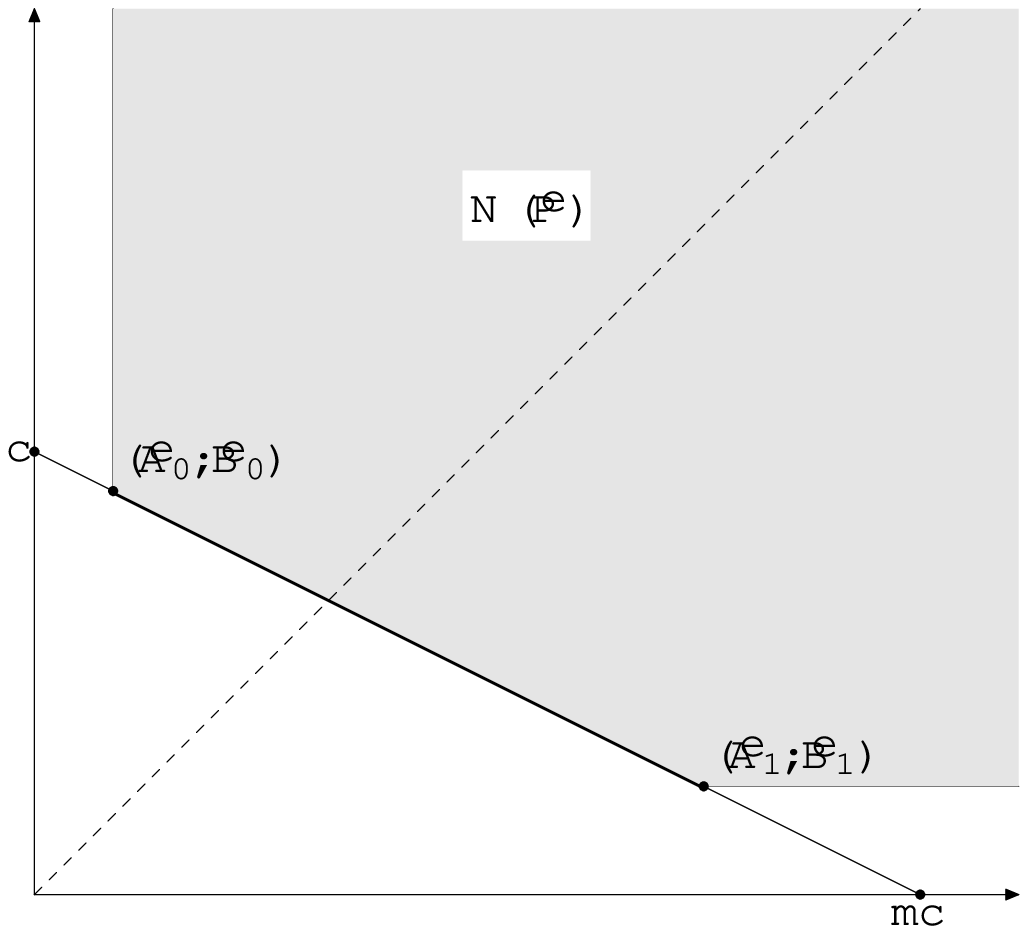}
\caption{\label{fig2}}
\end{figure}

\begin{lemma}\label{s3.2}
Let $m\in \bN,\, m\ge 1,$ be given, and denote by  $\mu$  the weight $\mu:=(1,m).$ Moreover, let $P$ be a $\mu$-homogeneous polynomial. Then its  Newton diagram $\N_d(P)$ is a compact  interval $[(A_0,B_0),(A_1,B_1)] $ joining two vertices  $(A_0, B_0), \, (A_1,B_1)$ (which may coincide), where we shall assume that $A_0\le A_1.$

Moreover, let $x=\vp(y)$ be a change of coordinates of the form $x_1=y_1$ and $x_2=y_2+a_2y_1^{m}$ (with $a_2\ne 0$),  if $m\ge 1,$ or $x_1=y_1+a_1y_2$ and $x_2=y_2+a_2y_1,$ if $m=1$ (with $a_j\ne 0$). Denote by $\tilde P$ the polynomial $P\circ\vp.$ Then $\tilde P$ is $\mu$-homogeneous of the same degree as $P,$ and its Newton diagram $\N_d(\tilde P)$ is an interval of the form $[(\tilde A_0,\tilde B_0),(\tilde A_1,\tilde B_1)],$ with $\tilde A_0\le \tilde A_1.$

\medskip
\noi  Assume either 
 
\bee 
\item[(i)] that the interval $\N_d(P)$ lies in the closed half-space above the bisectrix, i.e., $j\le k$ for every $(j,k)\in \N_d(P)$ (Figure \ref{fig1}), 

\item[(ii)]  or that  the point $(A_0, B_0)$ lies above or on the bisectrix, i.e., $A_0\le B_0,$ and  $m(P)\le d(P),$ where  $m(P)$ denotes again the maximal order of vanishing of $P$ along the unit circle $S^1,$ and $d(P)$ the distance. 
\ee

\noi Then the point $(\tilde A_0,\tilde B_0)$ lies in the closed half-space above the bisectrix and the point $(\tilde A_1,\tilde B_1)$ lies in the closed half-space below the bisectrix.
In particular, the Newton diagram $\N_d(\tilde P)$ intersects the bisectrix (Figure \ref{fig2})
 \end{lemma}

%%%%%%%%%%%%%%%%%%%%%%
%%%%%%%%%%%%%%%%%%%%%%
%%%%%%% Figure 1,2 %%%%%%%%%
%%%%%%%%%%%%%%%%%%%%%%
%%%%%%%%%%%%%%%%%%%%%%

 \noi {\bf Proof.} We have to show that the conditions (i) or (ii) imply
 \begin{equation}\label{n3.1}
 \tilde A_0\le \tilde B_0\ \mbox{and}\  \tilde A_1\ge \tilde B_1.
 \end{equation}
 
Now, by Proposition \ref{s2.1},  we can write $P$ in the form
\begin{equation}\label{3.8}
P(x_1,x_2)=x_1^{\al}x_2^{\beta}\prod_{l}(x_2-c_l
x_1^m)^{n_l},
\end{equation}
where the $c_l$'s are the non-trivial distinct complex roots of the polynomial $P(1, x_2)$ and the $n_l$'s  
are their multiplicities.

It is easy to read off  from \eqref{3.8} that the vertices of the  Newton polyhedron of $P$ are given by $(A_0,B_0)$ and $(A_1,B_1),$  
where 
 \begin{equation}\label{n3.2}
A_0:=\al,\, B_0:=\beta+N,\  \quad
A_1:=\al+mN,\quad B_1:=\beta,
\end{equation}
and that 
$$d_x=\frac{A_1+mB_1}{1+m}=\frac {\al+m(\beta+N)}{1+m}$$
(compare also \cite{PS} ); here, we have put
$$N:=\sum_{l}n_l.$$

\medskip
Assume first that $m\ge 1$  and $x_1=y_1,\, x_2=y_2+a_2y_1^{m}.$ Then, by \eqref{3.8}, we have 
 \begin{equation}\label{3.9}
\tilde P(y_1,y_2)=y_1^{\al}(y_2+a_2y_1^{m})^{\beta}\prod_{l}(y_2-(c_l-a_2)y_1^m)^{n_l}.
\end{equation}

 By looking at the terms of lowest  power of $y_1$ in \eqref{3.9}, we see that 
$\tilde A_0=\al,\, \tilde B_0=\beta+N,$ i.e., 
$$
(\tilde A_0, \tilde B_0)=( A_0, B_0),
$$
so that $\tilde A_0\le \tilde B_0.$ 
We next identify the  term of highest  power of $y_1$ in \eqref{3.9} in order to determine $(\tilde A_1,\tilde B_1)$.
\smallskip

If $a_2\ne c_l$ for every $l,$ then $\tilde A_1=\al+m\beta+mN,\  \tilde B_1=0$ hence $\tilde A_1>\tilde B_1.$

 Assume next that $a_2=c_{l_0}$ coincides with a non-trivial  root (which is then necessarily real).   Then we obtain
\begin{equation}\label{3.10}
\tilde A_1=\al+m\beta+m(N-n_{l_0}),\quad \tilde B_1=n_{l_0}.
\end{equation}
Now, if condition (i) holds, then by \eqref{n3.2} we have 
\begin{equation}\label{n3.3}
\al+mN\le \beta,
\end{equation}
which by Corollary \ref{s2.3} easily implies $n_{l_0}\le\al+m(\beta+N-n_{l_0}),$ and if (ii) holds, then 
$$n_{l_0}\le m(P)\le d_x=\frac {\al+m(\beta+N)}{1+m},$$ 
which is equivalent to  $n_{l_0}\le\al+m(\beta+N-n_{l_0}).$ Thus both conditions imply that 
\begin{equation}\label{3.11}
\tilde A_1\ge \tilde B_1,
\end{equation}
so that \eqref{n3.1} is satisfied. 

\medskip

There remains the case $ m=1$  and $x_1=y_1+a_1Y-2,\, x_2=y_2+a_2y_1^{m},$ with $a_1,a_2\ne 0.$ In this case,
\begin{equation}\label{3.12}
\tilde P(y_1,y_2)=(y_1+a_1y_2)^{\al}(y_2+a_2y_1)^{\beta}\prod_{l}\Big((1-c_la_1)y_2+(a_2-c_l)y_1\Big)^{n_l},
\end{equation}
where $a_1,a_2\ne 0.$  Our assumptions imply that either \eqref{n3.3} holds, or  that for every $l$ such that $c_l$ is real we have 
\begin{equation}\label{3.13}
n_l\le d_x=\frac{\al+\beta+N}2.
\end{equation}

Notice also that it may happen that $c_{l_1}a_1=1$ or $c_{l_0}=a_2,$ but at most for one $l_0$ and one $l_1,$ and these indices must be different, since $\vp$ has non-degenerate Jacobian at $0.$

Now, if $c_l a_1\ne 1$ for every $l,$ then 
$$
\tilde A_0=0, \ \tilde B_0=\al+\beta+N,
$$
and if $c_{l_1}a_1=1,$ then  
$$
\tilde A_0=n_{l_1}, \ \tilde B_0=\al+\beta+N-n_{l_1}.
$$
Moreover,  if $c_{l}\ne a_2$ for every $l,$ then 
$$
\tilde A_1=\al+\beta+N, \  \tilde B_1=0,
$$
and if $c_{l_0}=a_2,$ then  
$$
\tilde A_1=\al+\beta+N-n_{l_0}, \  \tilde B_1=n_{l_0}.
$$
Using  \eqref{n3.3} respectively \eqref{3.13}, it is easy to check that  $\tilde A_0\le \tilde B_0$ and  $\tilde A_1\ge \tilde B_1$ in all cases. 

\qed

 \bigskip
 Fix now a smooth function $f$ defined near the origin satisfying  $f(0,0)=0$ and
$\nabla f(0,0)=0.$

 Choose $\ka:=(\ka_1,\ka_2)$ with $\ka_j\ge 0$  such that  the principal face $\pi(f)$ of the Newton polyhedron lies on the line $\ka_1 t_1+\ka_2t_2=1.$  Notice that $\ka$ is uniquely determined, unless $\pi(f)$ is a vertex.  In the latter case, we shall assume that $\ka_1,\ka_2>0.$  Flipping coordinates $x_1$ and $x_2,$ if necessary, we shall assume without loss of generality that $\ka_1\le \ka_2.$ 
 Notice that then $\ka_1=0$ if and only if the principal face is non-compact; in this case, it is a half-line which lies on the horizontal line $L_\ka$ given by $t_2=1/\ka_2.$ And, if $\ka_1>0,$ then the  principal face is compact and the principal part $f_p=f_\ka$ of $f$ is a $\ka$-homogeneous polynomial of degree one. 

\medskip
Consider next a smooth local coordinate system $y$ at the origin given by  $x=\varphi(y),$
 and put 
 $$\tilde f(y_1,y_2):=f(\varphi(y_1,y_2)).$$
Flipping coordinates $y_1$ and $y_2,$ if necessary,we may  assume without loss of generality that 
$(x_1,x_2)=(\varphi_1(y_1,y_2), \varphi_2(y_1,y_2))$ satisfies 
$\frac{\pa \varphi_j(0,0)}{\pa y_j}\neq0$ for $j=1,2$. 

Therefore, we can write the functions   $\vp_1, \vp_2$ in the form
\begin{equation}\label{3.3}
\varphi_1(y_1,y_2)=y_1\psi_1(y_1,y_2)+\eta_1(y_2),\quad
\varphi_2(y_1,y_2)=y_2\psi_2(y_1,y_2)+\eta_2(y_1), 
\end{equation}

\noi where $\psi_1,\,\psi_2,\, \eta_1,\, \eta_2$ are smooth functions
satisfying
$$
\psi_1(0,0)\neq0,\quad \psi_2(0,0)\neq0, \quad
\eta_1(0)=\eta_2(0)=0.
$$
Since separate scaling of the coordinates $x_1$ and $x_2$ does not change the Newton polyhedron, 
we may further assume that $\psi_1(0,0)=\psi_2(0,0)=1.$ 

Denote by $k_j$ the order of vanishing of $\eta_j$ at $0,$ $j=1,2.$ Then clearly $k_j\ge 1.$  We shall see that if $k_2=\infty,$ i.e., if the function $\eta_2$ is flat at $0,$ then the distance will not change under the change of coordinates, i.e., $d_x=d_y.$ 

On the other hand, if $k_2$ is finite, then according to \eqref{3.3}, the main term of $\vp_2$ is of the form 
$y_2+a_2y_1^{k_2},$ with $a_2\ne 0.$ We therefore then  introduce a second  weight $\mu:=(1,k_2).$ Then the $\mu$-principal part   of $\vp_2(y_1,y_2)$  is given by $y_2+a_2y_1^{k_2},$ which is $\mu$-homogeneous of degree $k_2.$ Recall that  $L_\mu:=\{(t_1, t_2)\in\bR^2:t_1+k_2 t_2=d\}$ denotes the supporting line to the Newton polyhedron $\N(f).$  Notice that the line $L_\ka$ has slope $\frac {\ka_1}{\ka_2}\le 1,$ and  $L_\mu$ has slope $\frac 1{k_2}\le 1.$ The effect of the change of coordinates on the Newton diagram is then related to the interplay between these two homogeneities $\ka$ and $\mu,$  in particular the relation between the slopes of the corresponding supporting lines.

\medskip
We put $\frac{\ka_2}{\ka_1}:=\infty,$ if $\ka_1=0.$

\begin{lemma}\label{s3.3} 
In the situation described above, the following hold true:
\medskip

\noi a) Assume  either  that
\bee
\item[(i)]   $\ka_2>\ka_1$, and either $k_2=\infty,$  or $\infty>k_2>\frac{\ka_2}{\ka_1}$ (so that the line $L_\mu$ is less steep than the line $L_\ka$), or that 

\item[(ii)]  $\ka_1=\ka_2$ and $k_1, k_2>1.$
\ee
Then $d_y=d_x.$ 

\medskip

\noi a) If  $k_2<\frac{\ka_2}{\ka_1}$ (so that the line $L_\mu$ is steeper than the line $L_\ka$), then $d_y\le d_x.$ 

\medskip
In particular, if $\ka_1=0,$ i.e., if the principal edge is non-compact, then the coordinates $(x_1,x_2)$ are adapted to $f.$
\end{lemma}

\noi {\bf Proof.} {\bf a) The case where $\ka_1=0$ and $k_2=\infty.$ }
 \medskip
 
In this case $n:=\frac 1{\ka_2}\in\bN, \ n\ge 1, $ so that the  principal face is a half-line lying on the horizontal line $L_\ka$ given by $t_2=n,$ and $d_x=n$.  
Moreover,  $\eta_2$ is flat at the origin, i.e.,
\begin{equation}\label{3.4}
\partial_{ y_1}^j\eta_2(0)=0
\end{equation} 
for every $j\in\bN.$   In view of the structure of the Newton polyhedron that we assume, it is easy to see that then $f$ can be written in the form
$$f(x_1,x_2)=x_2^nx_1^{n_1}g(x_1)+x_2^{n+1}h(x_1,x_2)+\sum_{k=0}^{n-1}a_k(x_1)x_2^k,
$$
where $g$ and $h$ are smooth functions and $g(0)\ne 0,$ and where $a_1,\dots,a_{n-1}$ are flat functions at the origin.  Moreover, $n_1<n.$  This implies that 
$$\tilde f(y_1,y_2)= \vp_2(y_1,y_2)^n\, \vp_1(y_1,y_2)^{n_1}G(y_1,y_2)+\vp_2((y_1,y_2)^{n+1} H(y_1,y_2)+R(y_1,y_2)
$$
for some smooth function $G, H,$ and $R,$ where $G(0,0)\ne 0,$  and where  $R$ is flat at the origin. Since $\frac {\partial\vp_k}{\partial y_k}(0,0)=1,$  by  the product rule and \eqref{3.4} this easily implies that 
$\partial_{ y_1}^{n_1}\partial_{ y_2}^{n}\tilde f(0,0)\ne 0$ and $\partial_{ y_1}^{l}\partial_{ y_2}^{j}\tilde f(0,0)=0$ if $j<n,$ or if $ j=n$ and $l<n_1,$ so that $(n_1,n)\in \N(\tilde f),$  but no lattice point below or to the left of this point. This shows that  $d_y=d_x=n.$

\medskip
 \noi {\bf b) The case where $k_2>\frac{\ka_2}{\ka_1},$   $\ka_1>0,$ and either $\ka_2>\ka_1,$ or $\ka_1=\ka_2$ and $k_1>1.$ }
 \medskip
 
Let us first state some general observation. If $\vp_\ka$ denotes the $\ka$-principal part of $\vp,$ then it is easily seen, e.g.,  by Lemma \ref{s2.1n}, that 
$$
\tilde f(y_1,y_2)=f_\ka\circ\vp_\ka(y_1,y_2) +\ \mbox{terms of higher $\ka$-degree,}
$$
 so that
$$
\tilde f_\ka=f_\ka\circ\vp_\ka.
$$
 Moreover, $f_\ka\circ\vp_\ka$ is a $\ka$-homogeneous polynomial, so that its Newton diagram $\N_d(\tilde f_\ka)$ is again a compact interval (possibly a single point). In case that this interval intersects the bisectrix too, then it is  the principal face of $\N(\tilde f).$ Moreover, if $\vp_\ka$ has $\ka$-degree one, 
 then $\N_d(f_\ka)$ lies again on the line $L_\ka,$ and consequently we have $d_y=d_x.$ 
 
 Notice that the same conclusion holds true if $x_1,x_2$ are $\ka$-homogeneous of any degree $\delta,$ and if $\vp_\ka$ has the same $\ka$-degree $\delta.$
 
 \medskip
  Now, if  $k_2>\frac{\ka_2}{\ka_1},$ then $k_2\ka_1>\ka_2,$ and thus, by \eqref{3.3},  $x_2=y_2,$ up to terms of higher $\ka$-degree. Moreover, if $\ka_2>\ka_1$ or $k_1>1,$ then also $x_1=y_1,$ up to terms of higher $\ka$-degree, so that $\vp_\ka(y_1,y_2)=(y_1,y_2).$ Thus the reasoning above applies, and we see that $d_y=d_x.$

 \medskip
 \noi {\bf c) The case where $k_2<\frac{\ka_2}{\ka_1}$  and $\ka_1\ge 0.$}
\medskip

In this case,  the $\mu$-principal part   of $\vp_2(y_1,y_2)$  is given by $y_2+a_2y_1^{k_2},$ and the $\mu$-principal part of $\vp_1(y_1,y_2)$  is given by $y_1,$ if $k_2>1,$ and by $y_1+a_1y_2,$ if $k_2=1,$ with $a_1\ne 0$ if and only if   $k_1=1.$

Choose $d>0$ such that  the line $L_\mu:=\{(t_1, t_2)\in\bR^2:t_1+k_2 t_2=d\}$ is the supporting line to the Newton polyhedron $\N(f),$  and let $f_\mu$ be the $\mu$-principal part of $f.$   Then obviously $d\in\bN.$ Since the line $L_\mu$ is steeper than the line $\ka_1 t_1+\ka_2 t_2=1,$ it is clear from the geometry of the Newton polyhedron that $L_\mu \cap \N(f)$ is a compact interval $I_\mu$ (possibly a single point) lying in the closed half-space above the bisectrix, i.e., 

%\begin{equation}\label{3.6}
$$
f_\mu(x)=\sum_{(j,k)\in I_\mu}c_{jk}x_1^jx_2^k, \quad\mbox{where}\  j \le k\ \mbox{for every}\ (j,k)\in I_\mu.
$$
%\end{equation}
Applying Lemma \ref{s3.2} to the $\mu$-homogeneous polynomial $f_\mu,$ we see that the Newton diagram of $\tilde f_\mu=f_\mu\circ\vp_\mu$ intersects the bisectrix. The principal  face of $\N(\tilde f)$ lies therefore on the line $L_\mu,$ and since this line is steeper then the line $L_\ka$ containing the principal face of $\N(f),$ it is clear from the geometry that $d_y\le d_x.$ 

\medskip
Finally, assume that $\ka_1=0.$ Then either $k_2=\infty,$ or $k_2<\frac{\ka_2}{\ka_1},$ and so we always have $d_y\le d_x.$ Therefore,  the coordinates $(x_1,x_2)$ are adapted to $f.$

\qed

\medskip

\begin{thm}\label{s3.4}
Let $f$ be a real analytic (respectively smooth )  function near the origin, with $f(0,0)=0$ and
$\nabla f(0,0)=0. $ Assume that the given coordinates are not adapted  to $f.$ Then  all of the following conditions hold true:

\bee
\item[(a)] The principal face $\pi(f)$ of the Newton polyhedron is a compact edge. 

\noi It thus lies an a uniquely determined  line $\ka_1t_1+\ka_2 t_2=1, $ with $\ka_1,\ka_2 >0.$ 
Flipping coordinates $x_1,x_2,$ if necessary, we assume without loss of generality that $\ka_1\le \ka_2.$

\item[(b)]  $\frac{\ka_2}{\ka_1}\in\bN.$ 

\item[(c)] Any  root of $f_p$ of maximal order $m(f_p)$ on the unit circle $S^1$  lies away from the coordinate axes, and we have $m(f_p)>d(f).$  

\noi In particular, there exists a unique non-trivial real root $a$ of the polynomial $f_p(1,x_2)$ with 
multiplicity $n_a=m(f_p)>d(f).$
\ee
Moreover, in this   case, if we put $m:=\frac{\ka_2}{\ka_1}\in\bN,$ then  an adapted coordinate system for the principal part $f_p$ of $f$ is defined  by  
$y_1:=x_1,\,
y_2:=x_2-ax_1^m.$ The height  of $f_p$ is then given by $h(f_p)=m(f_p).$ 

Conversely, if all of the conditions (a) - (c) are satisfied, then the given coordinates are not adapted to the principal part $f_p$ of $f.$ 

\end{thm}

Notice that for mixed homogeneous polynomials $f=f_p$  the theorem gives a necessary and sufficient condition 
 for the adaptedness of the coordinates. We shall see later that the same is indeed true for general functions $f.$ 

Observe also that the root in (c) with multiplicity $n_a>d(f)$ is the principal root of $f_p.$

\medskip

\noi{\bf Proof of Theorem \ref{s3.4}.} 
Assume that the  coordinates $(x_1,x_2)$ are not adapted to $f.$ Then there exists a smooth  change of 
coordinates $x=\vp(y)$ at the origin such that 
\begin{equation}\label{3.6}
d_y>d_x.
\end{equation}
Arguing similarly as in the proof of Lemma \ref{s3.3}, we may assume that $\vp$ satisfies \eqref{3.3}, 
with 
$$
\psi_1(0,0)= \psi_2(0,0)=1, \quad
\eta_1(0)=\eta_2(0)=0.
$$
Denote again by $k_j$ the order of vanishing of $\eta_j$ at $0,$ $j=1,2,$ and choose the weight $\ka=(\ka_1,\ka_2)$ as before in such a way that the principal face of $\N(f)$ lies on the line $L_\ka.$ 
Again, assume that $\ka_1\le \ka_2.$ 

From \eqref{3.6} and Lemma \ref{s3.3} it then follows that if $\ka_2>\ka_1,$ then
$$\frac {\ka_2}{\ka_1}=k_2\in \bN,$$
and if $\ka_1=\ka_2,$ then $k_1=1$ or $k_2=1.$ In the latter case, by symmetry in the variables $x_1$ and $x_2,$ let us assume without loss of generality that $k_2=1.$ 
In particular,  the principal face $\pi(f)$ is compact, hence either a compact edge or a vertex. It is therefore an interval $[(A_0,B_0),(A_1,B_1)] $ joining two vertices  $(A_0, B_0), \, (A_1,B_1)$ (which may coincide), i.e., 
$$\pi(f)=[(A_0,B_0),(A_1,B_1)].$$
We shall assume that $A_0\le A_1.$ Since this interval  intersects the bisectrix $t_1=t_2,$ we then have 
\begin{equation}\label{3.7}
A_0\le B_0,\  A_1\ge B_1. 
\end{equation}

\medskip
In the sequel, let us write $m:=k_2,$  and  consider the $\ka$-homogeneous polynomial $f_p=f_\ka.$ Its Newton diagram is given by the interval $[(A_0,B_0),(A_1,B_1)] .$

\medskip
We next show that $m(f_\ka)>d_x.$
To this end, assume to the contrary that  $m(f_\ka)\le d_x,$ so that $f_\ka$ satisfies condition (ii) in Lemma \ref{s3.2}.
Observe  that the  $\ka$-principal part   of  $\vp_2(y_1,y_2)$  is given by $y_2+a_2y_1^{m},$ where $a_2\ne 0,$  and the $\ka$-principal part of $\vp_1(y_1,y_2)$  is given by $y_1,$ unless 
$\ka_1=\ka_2$  (hence $m=1$)  and $k_1=1,$ when it is given by $y_1+a_1y_2,$ with $a_1\ne 0.$ Therefore, the $\ka$-principal part $\tilde f_\ka$ of $\tilde f$ is given by
$$
\tilde f_\ka(y_1,y_2)=f_\ka(y_1+a_1y_2,y_2+a_2y_1^{m}),
$$
where $a_2\ne 0$ and $a_1\ne 0$ if and only if $m=k_1=1.$ Lemma \ref{s3.2} then shows that the Newton diagram of $\tilde f_\ka$ intersects the bisectrix, so that the principal face of $\N(\tilde f)$ lies on the line $L_\ka,$ hence $d_y=d_x.$ This contradicts \eqref{3.6}.

\medskip
We have thus seen that our assumptions on $f$ imply that $\frac {\ka_2}{\ka_1}\in \bN$ and $m(f_p)>d(f).$  
Assuming these properties, write the principal part $f_p=f_\ka$ of $f$ according to  Proposition \ref{s2.1} in the form 
\begin{equation}\label{3.8n}
f_\ka(x_1,x_2)=x_1^{\al}x_2^{\beta}\prod_{l}(x_2-c_l
x_1^m)^{n_l},
\end{equation}
where the $c_l$'s are the non-trivial distinct complex roots of the polynomial $f_\ka(1, x_2)$ and the $n_l$'s  
are their multiplicities.

Observe that clearly $\N(f_\ka)$ is contained in the half-plane $\{t_1\ge \al\}, $ so that $d_x=d(f_\ka)\ge \al.$ Similarly, we see that $d_x\ge \beta.$ 
Consequently,  there must be an  $l_0$ with real root $c_{l_0}$ so that $n_{l_0}=m(f_p).$  Notice that this root is unique, by Corollary \ref{s2.3}. Let us remark at this point that this  excludes the possibility that the principal face of $f$ is a vertex, for then
$f_p(x)$ would be of the form $cx_1^\al x_2^\beta,$ so that $f_p(1,x_2)$   had no non-trivial root.  Consequently, $\pi(f)$ is a compact edge.

\medskip
We show that in this case the change of coordinates $x_1:=y_1,\, x_2:=y_2+a_{l_0}y_1^m$ leads to adapted coordinates. To this end, we shall refer to the notation and results in the proof of Lemma \ref{s3.2}, applied to $P:=f_p.$

We have seen that  the function $\widetilde{ f_p}$ representing $f_p$ in the new coordinates $(y_1,y_2)$ is again a $\ka$-homogeneous polynomial, so that its Newton diagram is  a compact interval $[(\tilde A_0,\tilde B_0),(\tilde A_1,\tilde B_1)]$, whose endpoints are given by $(\tilde A_0, \tilde B_0)=( A_0, B_0)$ and \eqref{3.10}. Since  now $n_{l_0}> d_x,$ the order signs $\ge$ in \eqref{3.11} have to be replaced by $<,$ i.e., 
$$
\tilde A_1<\tilde B_1.
$$
This implies that the  interval $[(\tilde A_0,\tilde B_0),(\tilde A_1,\tilde B_1)]$ lies in the half-plane
$t_2>t_1.$ Consequently, the principal face of $\N(\widetilde{ f_p})$ lies on the horizontal line $t_2=\tilde B_1,$
hence  $d_y=\tilde B_1=m(f_p)$. Since this face is non-compact, the coordinates $(y_1,y_2)$ are adapted,  due to Lemma \ref{s3.3}. 

Notice, however, that the same change of coordinates will in general not lead to an adapted coordinate system for $f$ -- this may require further, higher order terms in $y_1,$ in addition to $a_{l_0}y_1^m.$

\medskip

Finally, assume conversely that all of the conditions (a) - (c) are satisfied. Then the coordinates $(x_1,x_2)$ are not adapted to $f_p,$ since, as we have seen, we can change coordinates for $f_p$ 
in such a way that, in the new coordinates $(y_1,y_2),$ the distance is given by $d_y=m(f_p)>d_x.$ 
\medskip
This concludes the proof of Theorem \ref{s3.4}.

\qed

As a corollary, we obtain the following characterization of the height in case of a $\ka$-homogeneous polynomial.
\begin{cor}\label{s3.5}
Let $P$ be a $\ka$-homogeneous polynomial as in Proposition \ref {s2.1}. Then
$$h(P)=\max\{m(P),d_h(P)\}.
$$
\end{cor}

\proof We adopt the notation from Proposition \ref{s2.1}. If $\frac{\ka_2}{\ka_1}\notin \bN,$ then the coordinates are adapted to $P,$ and  by Corollary \ref{2.3} (a) and \eqref{2.9} 
we have $\max\{m(P),d_h(P)\}=\max\{\nu_1,\nu_2,d_h(P)\}=d(P)=h(P).$

So, assume next that $\frac{\ka_2}{\ka_1}\in \bN.$ Then, according to Corollary \ref{2.3} (b), there is at most one real root of $P$ on the unit circle with multiplicity $n_{l_0}>d_h(P).$ 

If there is no such root, then $m(P)\le d_h(P),$ and so the coordinates are adapted. This shows that $h(P)=d(P)=d_h(P)=\max\{m(P),d_h(P)\}.$

If there is such a root, and if it lies on one of the coordinate axes, then $n_{l_0}=\nu_1$ or $n_{l_0}=\nu_2,$ hence $m(P)=\max\{\nu_1,\nu_2\},$ and the claim follows. We may thus assume that the root does not lie on a coordinate axis, so that $1\le l_0\le M.$ Then $\nu_1,\nu_2\le d_h(P),$ which implies that the principal face $\pi(P)$ of the Newton polyhedron of $P$ must be a compact edge. Then, by Theorem \ref{s3.4}, the coordinates are not adapted to $P,$ and  $h(P)=m(P)>d(P).$ Since $d(P)\ge d_h(P),$ the conclusion follows also in this case.

\qed

\bigskip

\setcounter{equation}{0}
\section{The real analytic case}\label{analytic}

In this section we prove that if $f$ is a real analytic function  then there exists an 
analytic  coordinate system which is adapted to $f.$  This coordinate
system can be described  in terms  of the  Puiseux series expansion of the roots of the
equation $f\x=0.$

\medskip
\subsection{Description of the Newton polyhedron in terms of the roots.}

Assume again  $f$ to be real analytic and real valued. By  the Weierstra\ss{}  preparation theorem  we can then write
$$
f(x_1,x_2)=U(x_1,x_2)x_1^{\nu_1} x_2^{\nu_2} F(x_1,x_2)
$$
near the origin, where $F(x_1,x_2)$ is a pseudopolynomial of the form
$$
F(x_1,x_2)=x_2^m+g_1(x_1)x_2^{m-1}+\dots+ g_m(x_1),
$$
and $U,\, g_1,\dots , g_m$ are real analytic functions satisfying 
 $U(0,0)\neq0, \, g_j(0)=0$. Observe that the Newton polyhedron of $f$ is the same as that of $x_1^{\nu_1} x_2^{\nu_2} F(x_1,x_2).$ We shall also assume without loss of generality that $g_m$ is a non-trivial function, so that the roots $r(x_1)$ of the  equation $F(x_1,x_2)=0,$ considered as a polynomial in $x_2,$ are all non-trivial. 
 
Following \cite{PS}, it is well-known  that these roots can be expressed in a small neighborhood of $0$ as   Puiseux series 
$$
r(x_1)=c_{l_1}^{\al_1}x_1^{a_{l_1}}+c_{l_1l_2}^{\al_1\al_2}x_1^{a_{l_1l_2}^{\al_1}}+\cdots+
c_{l_1\cdots l_p}^{\al_1\cdots \al_p}x_1^{a_{l_1\cdots
l_p}^{\al_1\cdots \al_{p-1}}}+\cdots,
$$
where
$$
c_{l_1\cdots l_p}^{\al_1\cdots \al_{p-1}\be}\neq c_{l_1\cdots
l_p}^{\al_1\cdots \al_{p-1}\ga} \quad \mbox{for}\quad \be\neq \ga,
$$

$$
a_{l_1\cdots l_p}^{\al_1\cdots \al_{p-1}}>a_{l_1\cdots
l_{p-1}}^{\al_1\cdots \al_{p-2}},
$$
with strictly positive exponents $a_{l_1\cdots l_p}^{\al_1\cdots \al_{p-1}}>0$ and non-zero complex coefficients $c_{l_1\cdots l_p}^{\al_1\cdots \al_p}\ne 0,$ and where we have kept enough terms to distinguish between all the non-identical roots of $F.$

By the {\it cluster } $\left[ \begin{matrix} \al_1&\cdots &\al_p\\
l_1&\dots &l_{p}
\end{matrix}\right ],$
we shall designate all the roots $r(x_1)$, counted with  their multiplicities,
which satisfy
\begin{equation}\label{4.1}
r(x_1)-c_{l_1}^{\al_1}x_1^{a_{l_1}}+c_{l_1l_2}^{\al_1\al_2}x_1^{a_{l_1l_2}^{\al_1}}+\cdots+
c_{l_1\cdots l_p}^{\al_1\cdots \al_p}x_1^{a_{l_1\cdots
l_p}^{\al_1\cdots \al_{p-1}}}=O(x_1^b) \
\end{equation}
for some exponent $b>a_{l_1\cdots l_p}^{\al_1\cdots \al_{p-1}}$. We
also introduce the clusters 
 $\left[ \begin{matrix}
  \al_1&\cdots &\al_{p-1}&\cdot\\
l_1&\dots &l_{p-1}& l_p
\end{matrix}\right ],$
by
$$
 \left[ \begin{matrix} 
 \al_1&\cdots &\al_{p-1}&\cdot\\
l_1&\dots &l_{p-1}& l_p
\end{matrix}\right ]
:=
\bigcup\limits_{\al_p}
\left[ \begin{matrix} 
\al_1&\cdots &\al_p\\
l_1&\dots &l_p
\end{matrix}\right ].
$$

Each index $\al_p$ or $l_p$ varies in some finite range which we shall not specify here. We finally put 
$$
N\left[ \begin{matrix} 
\al_1&\cdots &\al_p\\
l_1&\dots &l_p
\end{matrix}\right ]:=\mbox{number of roots in} \,
\left[ \begin{matrix} 
\al_1&\cdots &\al_p\\
l_1&\dots &l_p
\end{matrix}\right ],
$$
$$
N \left[ \begin{matrix} 
 \al_1&\cdots &\al_{p-1}&\cdot\\
l_1&\dots &l_{p-1}& l_p
\end{matrix}\right ]:=\mbox{number of roots in}
\,  \left[ \begin{matrix} 
 \al_1&\cdots &\al_{p-1}&\cdot\\
l_1&\dots &l_{p-1}& l_p
\end{matrix}\right ]
$$
\medskip

Let $a_1<\dots< a_l<\dots<a_n$ be the distinct leading exponents of all the 
roots of $F.$ Each exponent $a_l$ corresponds to the cluster $\left[ \begin{matrix} 
\cdot\\
l
\end{matrix}\right ],$ so that  the set of all roots of $F$ can be divided as $\bigcup\limits_{l=1}^n
\left[ \begin{matrix} 
\cdot\\
l
\end{matrix}\right ]$.
Then we may write
$$
f(x_1,x_2)=U(x_1,x_2)x_1^{\nu_1} x_2^{\nu_2} \prod_{l=1}^n
\Phi\left[ \begin{matrix} 
\cdot\\
l
\end{matrix}\right ](x_1,x_2),
$$
where
$$
\Phi\left[ \begin{matrix} 
\cdot\\
l
\end{matrix}\right ](x_1,x_2):=\prod_{r\in \left[ \begin{matrix} 
\cdot\\
l
\end{matrix}\right ]}(x_2-r(x_1)).
$$

We introduce the following quantities:
\begin{equation}\label{4.2}
A_l=A\left[ \begin{matrix} 
\cdot\\
l
\end{matrix}\right ]:=\nu_1+\sum_{\mu\le l}a_\mu N\left[ \begin{matrix} 
\cdot\\
\mu
\end{matrix}\right ],\quad
B_l=B\left[ \begin{matrix} 
\cdot\\
l
\end{matrix}\right ]:=\nu_2+\sum_{\mu\ge l+1}N\left[ \begin{matrix} 
\cdot\\
\mu
\end{matrix}\right ].
\end{equation}
Notice that $B_l$ is just the number of all roots with leading exponent strictly  greater than $a_l$ (where we here interpret the trivial roots corresponding to the factor $(x_2-0)^{\nu_2}$ in our representation of $f\x$ as roots with exponent $+\infty$). 

If $N_{\all}$ denotes the total number of all roots of $f$ away from the axis $x_1=0,$ including the trivial ones and  counted with their multiplicities, then  we can also write 
\begin{equation}\label{4.3}
B_l=N_\all-\sum_{\mu\le l}N\left[ \begin{matrix} 
\cdot\\
\mu
\end{matrix}\right ].
\end{equation}

Then, similarly as for the reduced Newton diagram in \cite{PS}, the vertices of the Newton diagram $\N_d(f)$ of $f$ are the  points
$(A_l,B_l), \  l=0,\dots,n,$ and the Newton polyhedron $\N(f)$ is the convex hull
of the set $\cup_{l} ((A_l,B_l)+\bR_+^2)$.

Notice also that 
$$A_l+a_lB_l=A_{l-1}+a_lB_{l-1}.
$$

Let $L_l:=\{(t_1,t_2)\in \bN^2:\ka^l_1t_1+\ka^l_2 t_2=1\}$ denote the line passing through the  points $(A_{l-1},B_{l-1})$ and $(A_l,B_l).$ It is easy to see that it is given by
\begin{eqnarray}\nonumber
\ka^l_1&=&\frac {\Delta B_l}{A_l\Delta B_l-B_l\Delta A_l}=\frac 1{A_l+a_lB_l}\nonumber \\
&&\\ \label{4.4}
\ka^l_2&=&\frac {\Delta A_l}{A_l\Delta B_l-B_l\Delta A_l}=\frac {a_l}{A_l+a_lB_l},\nonumber
\end{eqnarray}
where $\Delta A_l:=A_l-A_{l-1}, \Delta B_l:=B_l-B_{l-1}.$  
This implies that 
\begin{equation}\label{4.5}
\frac {\ka^l_2}{\ka^l_1} =a_l,
\end{equation}
which in return is the reciprocal of the slope of the line $L_l.$  The line $L_l$ intersects the bisectrix at the point  $(d_l,d_l), $ where 
$$d_l :=\frac{A_l+a_lB_l}{1+a_l}.$$ 
Moreover, the vertical edge of the Newton polyhedron, which passes through the point 
$(A_0,B_0)=(\nu_1,\nu_2 +m), $ intersects the bisectrix at $(\nu_1,\nu_1),$ and the horizontal edge, which passes through the point $(A_n,B_n)=(A_n, \nu_2),$ intersects the bisectrix at $(\nu_2,\nu_2).$ We therefore conclude that the distance $d(f)$  is given by 
$d(f)=\max \{\nu_1, \nu_2, \max_{l =1, \dots, n} d_l\}.$ 

Finally, fix $l,$ and let us determine the $\ka^l$-principal part $f_{\ka^l }$ of $f$ corresponding to the supporting line $L_l.$ To this end, observe that $f$ has the same $\ka^l$-principal part  as the function 
$$
U(0,0)x_1^{\nu_1}x_2^{\nu_2}\prod_{\al,\mu}\Big(x_2-c^\al_\mu x_1^{a_{\mu}}\Big)^{N \left[ \begin{matrix} 
\al\\
\mu
\end{matrix}\right ]}.
$$
Moreover, the $\ka^l$-principal part of $x_2-c^\al_\mu x_1^{a_{\mu}}$ is given by $c^\al_\mu x_1^{a_{\mu}},$ if $\mu<l, $ and by $x_2,$ if $\mu>l.$ This implies that 
\begin{equation}\label{4.5}
f_{\ka^l }\x=c_l \,x_1^{A_{l-1}} x_2^{B_l}\prod_\al \Big(x_2-c^\al_l x_1^{a_{l}}\Big)^{N\aol}.
\end{equation}
In view of this identity, we shall say that the edge $[(A_{l-1},B_{l-1}) ,(A_l,B_l)]$ is {\it associated}  to the cluster of roots $\dotol.$
We collect these results in the following lemma.

\begin{lemma} \label{s4.1} 
The vertices of the Newton diagram $\N_d(f)$ of $f$ are the  points
$(A_l,B_l), \  l=0,\dots,n,$ with $A_j, B_j$ given by \eqref{4.2}, and its  edges are the intervals 
$[(A_{l-1},B_{l-1}) ,(A_l,B_l)],$ $l=1,\dots,n.$  Moreover, the distance between  the Newton polyhedron and the origin is given by 
\begin{equation}\label{4.6}
d(f)=\max \left\{A_0,B_n,\max_{l =1, \dots, n}\frac{A_l+a_lB_l}{1+a_l}\right\},
\end{equation}
and the $\ka^l$-principal part  of $f$ corresponding to the supporting line $L_l$ through the edge $[(A_{l-1},B_{l-1}) ,(A_l,B_l)]$ is given by \eqref{4.5}.
\end{lemma}

\medskip
\subsection{Existence of adapted coordinates in the analytic setting.}

\begin{thm} \label{s4.2}
Let $f$ be a non-trivial real valued  real analytic function defined on a neighborhood of the origin in $\bR^2$ with
$f(0,0)=0,\,\nabla f(0,0)=0$. Choose $\ka_1,\ka_2\ge 0$ such that the principal face $\pi(f)$ of
the Newton polyhedron  of $f$  lies on the   line $\ka_1 t_1+\ka_2
t_2=1.$  Without loss of generality,  we may assume that  $\ka_2\ge\ka_1$. Then there exists a  real analytic function 
$\psi(x_1)$ of $x_1$ near the origin with $\psi(0)=0$ such that an adapted coordinate
system $(y_1,y_2)$ for $f$ near $0$ is given by  $y_1:=x_1,\, y_2:=x_2-\psi(x_1).$ 

\smallskip

The function $\psi$ can in fact be chosen as one of the roots, respectively a leading partial sum of the Puiseux series expansion of one of the roots, of the equation $f\x=0,$ considered as an equation in $x_2.$
\end{thm}

Notice  that  Theorem \ref{s4.1} was proved by A.N. Varchenko in \cite{Va}.
But, his proof is based on H.~Hironaka's  \cite{H} deep theorem on the 
resolution of singularities. We shall give a more elementary proof  of
Varchenko's theorem,  based on the Puiseux series expansion of roots of $f,$  which in fact gives 
an explicit description of an  adapted coordinate system in terms of these roots.

\medskip
\noi{\bf Proof.}  We may and shall assume without loss of generality that $U\equiv 1.$ If the coordinate system is adapted, then we can choose $\psi=0.$  According to Theorem \ref{s3.4}, this applies in any of the following three cases:
\bee
\item[(a)] The principal face is unbounded. Since we are assuming that $\ka_2\ge \ka_1,$ by Lemma \ref{s4.1} this means that $A_n<B_n,$ i.e., that 
\begin{equation}\label{4.8}
\nu_1+\sum_{l=1}^n a_l N\dotol <\nu_2.
\end{equation}

\item[(b)] $\pi(f)$ consists of  a vertex. Then choose $0\le \la\le n$ so that $\pi(f)=\{(A_\la,A_\la)\}.$  This happens iff 
\begin{equation}\label{4.9}
\nu_1+\sum_{l\le \la} a_l N\dotol =\nu_2+\sum_{l\ge \la+1}  N\dotol.
\end{equation}

\item[(c)] $\pi(f)$ is a compact edge $[(A_{\la-1},B_{\la-1}) ,(A_\la,B_\la)], $ but 
\begin{equation}\label{4.10}
a_\la=\frac {\ka^{\la}_2}{\ka^\la_1} \notin\bN
\end{equation}
or 
\begin{equation}\label{4.11}
a_\la\in \bN \ \mbox{and}\ 
N\left[ \begin{matrix} 
\al\\
\la
\end{matrix}\right ]\le d(f)\quad\mbox{whenever}\quad c^{\al}_\la\in\bR.
\end{equation}

\begin{figure}
\centering
\includegraphics[scale=0.7]{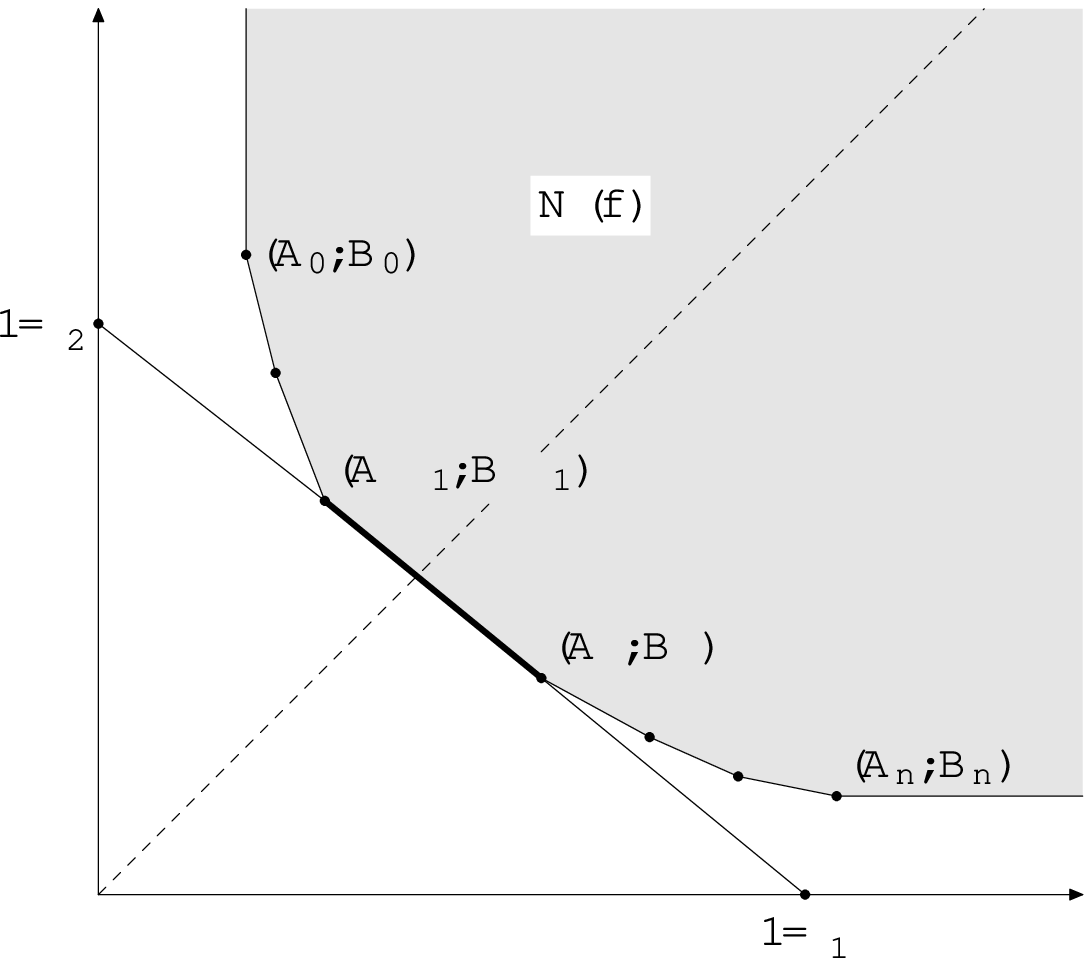}
\caption{\label{fig3}}
\end{figure}

\ee
\medskip
There remains the case where  the principal face $\pi(f)$ is a  compact edge
$\pi(f)= [(A_{\la-1},B_{\la-1}) ,(A_\la,B_\la)]$ ( $1\le \la\le n$), $a_\la\in\bN$
and there is  an index $\beta$ such that 
\begin{equation}\label{4.12}
m(f_p)=N\bola>d(f)=\frac{A_{\la}+a_{\la} B_{\la}}{1+a_{\la}}\quad\mbox{and}\quad c^{\beta}_\la\in\bR
\end{equation} 
(Figure \ref{fig3}). Notice that in view of Corollary \ref{s2.3} the index $\beta$ is unique, and $c_\la^\beta x_1^{a_\la}$ is the principal root of $f_p.$
\smallskip

We then apply an algorithm due to  Varchenko \cite{Va}.
\medskip

%%%%%%%%%%%%%%%%%%%%%%%%%%%%%
%%%%%%%%%%%%%%%%%%%%%%%%%%%%%
%%%%%%%% Figure3,4 %%%%%%%%%%%%%%%%
%%%%%%%%%%%%%%%%%%%%%%%%%%%%%
%%%%%%%%%%%%%%%%%%%%%%%%%%%%%

\begin{figure}
\centering
\includegraphics[scale=0.7]{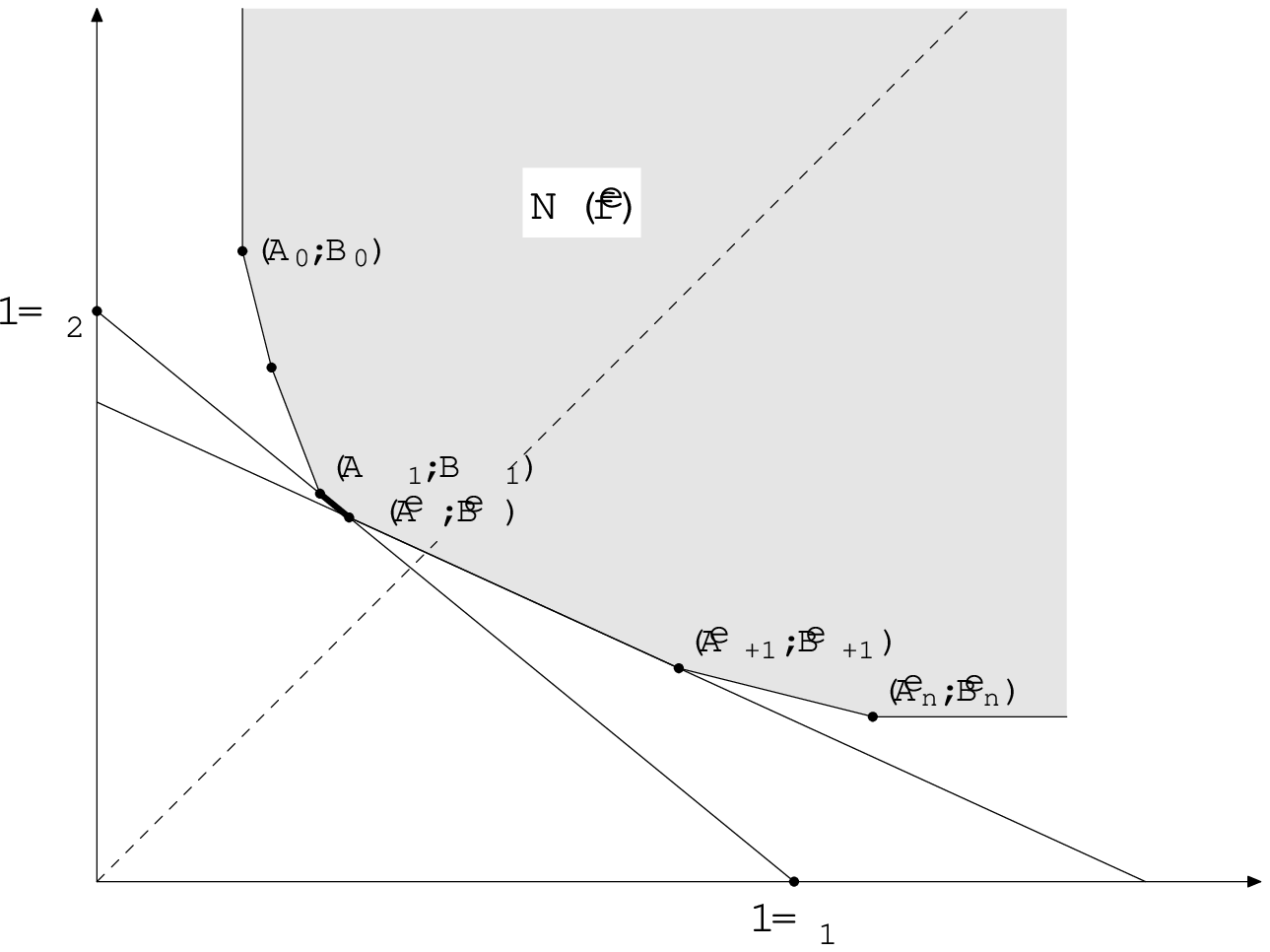}
\caption{\label{fig4}}
\end{figure}

\noi{\bf Step 1 (Figure \ref{fig4}).} We apply the  real  change of variables $x=\vp(y)$ given by $y_1:=x_1,\, y_2:= x_2-c_\la^\beta x_1^{a_\la},$  and put again $\tilde f:=f\circ \vp.$ In order to describe the effect of this change of variables to the Newton polyhedron, let us denote all quantities associated to $\tilde f$ with a superscript $\widetilde{}\,$.

Observe first that the roots $\tilde r$ of $\tilde f$ are of the form
\begin{equation}\label{4.13}
\tilde r(y_1)=c_{l}^{\al_1}y_1^{a_{l}}-c_\la^\beta
y_1^{a_{\la}}+c_{l l_2}^{\al_1\al_2}y_1^{a_{ll_2}^{\al_1}}+\cdots,
\end{equation}
respectively $\tilde r(y_1)=-c_\la^\beta y_1^{a_\la} $ (namely those corresponding to the trivial roots of $f$).  This shows that the leading exponents 
$$\tilde a_1<\tilde a_2< \cdots < \tilde a_l<\cdots
$$
of these roots are given  by $\tilde a_l=a_l$ for $l<\la,$ with the same multiplicities $N\widetilde\dotol=N\dotol$ (where $\widetilde\dotol$ denotes the cluster of roots $\tilde r$ associated to the index $l$).

For the vertices $(\tilde A_l,\tilde B_l)$ of $\N(\tilde f)$ we thus obtain from  \eqref{4.2}, \eqref{4.3} that 
\begin{equation}\label{4.14}
(\tilde A_l,\tilde B_l)=(A_l,B_l)\quad\mbox{for}\quad l<\la.
\end{equation}
Moreover, any root $r$ belonging to a cluster $\dotol$ with $l>\la$ is transformed into a root $\tilde r$ 
with leading exponent $a_\la.$ Finally, if $r$ belongs to $\dotola,$ then either the leading exponent of $\tilde r$ is $a_\la$ (namely if $\al_1\ne \beta$), or it is of the form $a_{\la l_2}^{\beta}>a_\la.$ 
We therefore distinguish two cases.

\medskip

\noi {\bf Case 1.} $\nu_2+B_\la+\Big(N\dotola-N\bola\Big)>0.$  

\medskip
\noi Then there exists at least one root $\tilde r$ with leading exponent $a_\la,$ so that $\tilde a_\la=a_\la.$ Moreover, since  $\tilde B_\la$ is the number of roots $\tilde r$ with leading exponent strictly greater than $a_\la,$ we  see that $\tilde B_\la=N\bola.$ 
Similarly, the number $N\widetilde\dotola$ is the same as the number of roots $ r$ with leading exponent strictly greater than $a_\la$ or equal to $a_\la,$ but then with index $\al_1\ne \beta,$  hence $N\widetilde\dotola=B_{\la-1}-N\bola.$ This implies $\tilde A_\la=\tilde A_{\la-1}+\tilde a_\la N\widetilde\dotola=A_{\la-1}+a_\la\Big(B_{\la-1}-N\bola\Big)=A_\la +a_\la B_\la-a_\la N\bola.$ In combination, we thus have 
\begin{equation}\label{4.14n}
(\tilde A_\la,\tilde B_\la)=\Big(A_\la +a_\la B_\la-a_\la N\bola,N\bola\Big).
\end{equation}

But, estimate \eqref{4.12} is equivalent to $\tilde A_\la<\tilde B_\la,$ so that the edge $[(\tilde A_{\la-1},\tilde B_{\la-1}),(\tilde A_\la,\tilde B_\la)],$ which lies on the same line as the edge $\pi(f)= [(A_{\la-1},B_{\la-1}) ,(A_\la,B_\la)]$ and has the same left vertex, is not the principal face of the Newton polyhedron of $\tilde f.$  Finally, it is evident from \eqref{4.13} that in this case
$\tilde a_{\la+k}=a_{\la k}^\beta$ if $ k>0$ (unless there is no cluster $\dotol$ with  $l>\la$)

This shows that in this case the principal face of $\N(\tilde f)$ is either associated to  a cluster of roots $\tilde r$ which corresponds to a cluster of roots  $\left[ \begin{matrix} \beta&\cdot\\
\la&\la_2
\end{matrix}\right ]$  in the original coordinates, or is a horizontal, unbounded edge (so that the new coordinates are adapted). 

\medskip 

\noi {\bf Case 2.} $\nu_2+B_\la+\Big(N\dotola-N\bola\Big)=0.$  

\medskip
\noi Then there is no root $\tilde r$ with leading exponent $a_\la,$ and so again the conclusion stated at the end of the previous case applies. 

In both cases we see that the principal face of the new Newton polyhedron $\N(\tilde f)$ will be less steep than the one of $\N(f),$ so that $d(\tilde f)>d(f).$ 
\medskip

\noi{\bf Subsequent steps.} 
Now, either the new coordinates $y$ are adapted, in which case we are finished. Or we can apply the same procedure to $\tilde f.$  Composing the change of coordinates from the first step with the one from the second step, we see that we then can find a change of coordinates $x=\vp_{(2)}(y)$ of the form
$$x_1:=y_1, \, x_2:= y_2-(c_\la^\beta x_1^{a_\la}+c_{\la\la_2}^{\beta \beta_2}y_1^{a_{\la\la_2}^\beta}),
$$ 
 with $a_{\la\la_2}^{\beta}\in\bN$ and $c_{\la\la_2}^{\beta \beta_2}\in \bR, $ such that the following holds:
\medskip
If the function $f_{(2)}:=f\circ\vp_{(2)}$  expresses the function $f$ in the  new coordinates, then the principal face of the Newton polyhedron of $f_{(2)}$ is either associated to  a cluster of roots which corresponds to a cluster of roots  $\left[ \begin{matrix} \beta&\beta_2&\cdot\\
\la&\la_2&\la_3
\end{matrix}\right ]$  in the original coordinates, or is a horizontal, unbounded edge (so that the new coordinates are adapted). 

\medskip
Now, if  we iterate this procedure, then either  this procedure will stop after finitely many steps, or it will continue infinitely. If it stops, it is clear that we will have arrived at a new, adapted coordinate system of the form
$$
x_1:=y_1, \, x_2:= y_2-(c_\la^\beta y_1^{a_\la}+\dots +c_{\la\la_2\cdots \la_p}^{\beta \beta_2\cdots\beta_p}y_1^{a_{\la\la_2\cdots\la_p}^{\beta\beta_1\cdots\beta_{p-1}}}),
$$ 
and Theorem \ref{s4.2} holds with a polynomial function $\psi(y_1).$  

\medskip

\begin{figure}[t]
\centering
\includegraphics[scale=0.7]{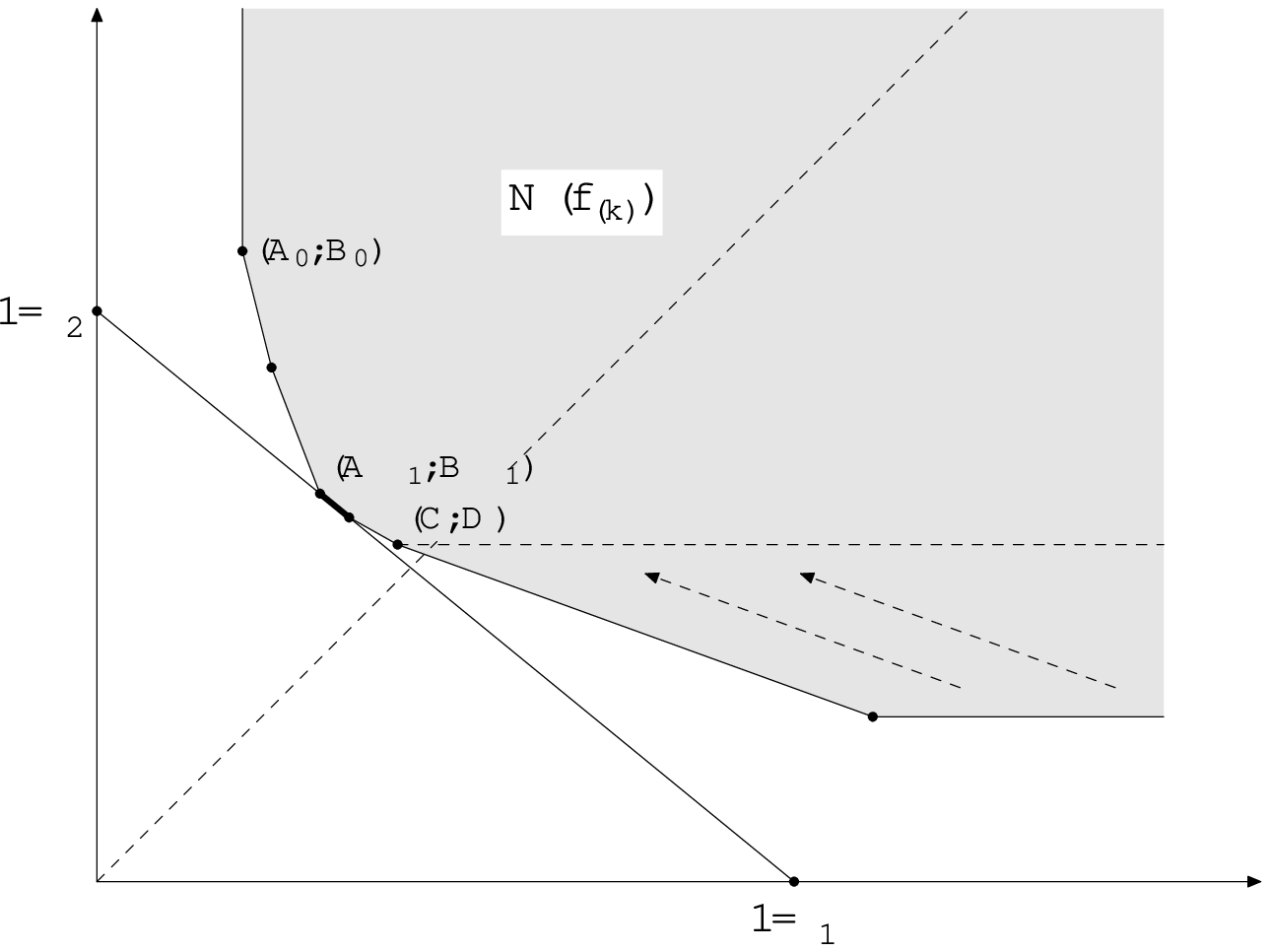}
\caption{\label{fig5}}
\end{figure}

\begin{figure}[h]
\centering
\includegraphics[scale=0.7]{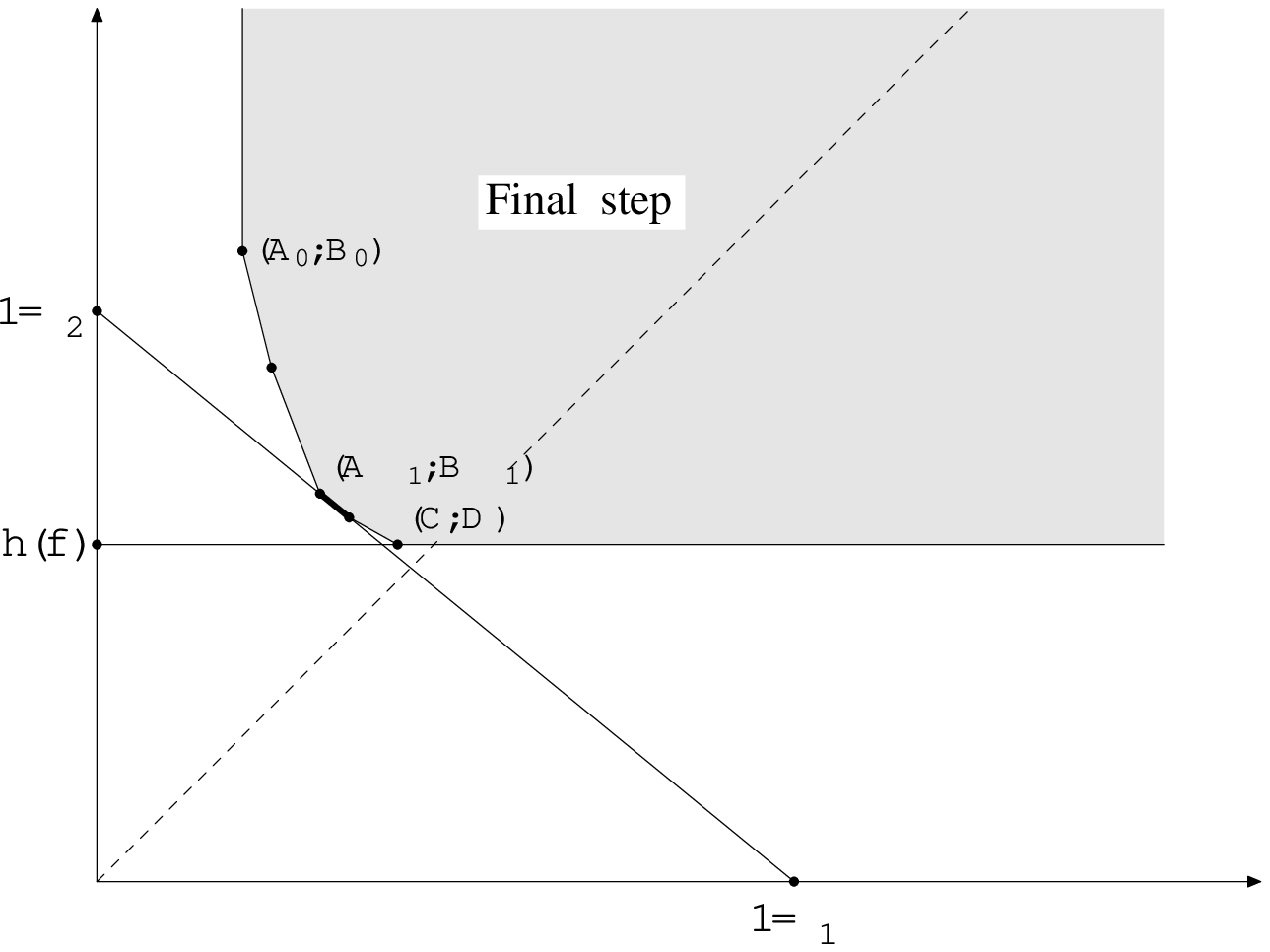}
\caption{\label{fig6}}
\end{figure}

\noi{\bf Final step (Figures \ref{fig5},\ref{fig6}).} Assume that the procedure does not terminate. Then, in each further step we pass to a new Newton diagram with a  principal face associated  to a cluster of roots (in the old coordinates) 
$\left[ \begin{matrix}\beta &\beta_2&\cdots &\cdot\\
\la&\la_2&\dots &\la_{k+1}
\end{matrix}\right ],$ which is a sub-cluster of the previous cluster $\left[ \begin{matrix}\beta &\beta_2&\cdots &\cdot\\
\la&\la_2&\dots &\la_{k}
\end{matrix}\right ].$ In particular, the corresponding muliplicities $N_k:=N\left[ \begin{matrix}\beta &\cdots &\cdot\\
\la&\dots &\la_{k+1}
\end{matrix}\right ]$
form a decreasing sequence, so that they eventually must become constant. Choose $N,k_0\in\bN$ such that $N_k=N$ for every $k\ge k_0.$ Replacing the original coordinate system by the one obtained in the $k_0$-th step, we may assume without loss of generality that $k_0=0,$ i.e., 
\begin{equation}\label{4.15}
N=N\dotola=N\left[ \begin{matrix} \beta&\cdot\\
\la&\la_2
\end{matrix}\right ]
=\dots=N\left[ \begin{matrix}\beta &\beta_2 &\cdots &\cdot\\
\la&\la_2&\dots &\la_{k+1}
\end{matrix}\right ]=\cdots
\end{equation}
This clearly implies that 
$N\left[ \begin{matrix}\beta &\beta_2 &\cdots &\cdot\\
\la&\la_2&\dots &\la_{k+1}\end{matrix}\right ]=N\left[ \begin{matrix}\beta &\beta_2 &\cdots &\beta_{k+1}\\
\la&\la_2&\dots &\la_{k+1}\end{matrix}\right ]
$
for every $k,$  and that each of the clusters $\left[ \begin{matrix}\beta &\beta_2 &\cdots &\cdot\\
\la&\la_2&\dots &\la_{k+1}\end{matrix}\right ]$ contains exactly  one and the same root (of multiplicity $N$), namely
\begin{equation}\label{4.16}
\rho(x_1):=c_\la^\beta x_1^{a_\la}+\dots +c_{\la\la_2\cdots \la_{k+1}}^{\beta \beta_2\cdots\beta_{k+1}}x_1^{a_{\la\la_2\cdots\la_{k+1}}^{\beta\beta_1\cdots\beta_{k}}}+\cdots .
\end{equation}
(so that in fact $\la_k=\beta_k=1$ for every $k\ge 2.$)
Moreover, our procedure shows that all cofficients in this series must be real, and all exponents positive integers, so that $\rho(x_1)$ is a real valued, real analytic function of $x_1.$  

%%%%%%%%%%%%%%%%%%%%%%%%%
%%%%%%%%%%%%%%%%%%%%%%%%%
%%%%%% Figures 5,6 %%%%%%%%%%%%%
%%%%%%%%%%%%%%%%%%%%%%%%%
%%%%%%%%%%%%%%%%%%%%%%%%%

If we apply  our first change of coordinates $y_1:=x_1, \, y_2:= x_2-c_\la^\beta x_1^{a_\la}$ in this situation, then we see that the leading exponents of the new roots $\tilde r$ are given by $ \tilde a_1=a_1,\dots, \tilde a_\la=a_\la$ and $\tilde a_{\la+1}=a_{\la \la_2}^\beta$ in  Case 1, and  by $ \tilde a_1=a_1,\dots, \tilde a_{\la-1}=a_{\la-1}$ and $\tilde a_{\la}=a_{\la \la_2}^\beta$  in Case 2. 
Moreover, it is clear from our discussion that the last edge, associated to the cluster of roots with biggest leading exponent $a_{\la \la_2}^\beta,$ must be  the principal edge. Moreover, in the new coordinates $y,$ the function $\tilde  f$ can have no vanishing  root, since such a root would have corresponded to a root $c_\la^\beta x_1^{a_\la},$ which, by \eqref{4.16}, cannnot exist in the cluster $\dotola.$ 

\medskip 

By passing to these new coordinates, we may therefore assume in addition to \eqref{4.15} that 
the cluster $\dotola$ is the cluster associated to the principal face of the Newton polyhedron  of $f,$ and that $\nu_2=0,$ i.e., that 
\begin{equation}\label{4.17a}
f_p\x=cx_1^{\nu_1} (x_2-c^\beta_\la x_1^{a_\la})^N.
\end{equation}
 The Newton polyhedron of $f$ thus has vertices $(A_0,B_0),\dots,(A_{\la},B_{\la}),$ and the principal edge is given by $[(A_{\la-1},B_{\la-1}) ,(A_\la,B_\la)],$ where, according to \eqref{4.2},  \eqref{4.15},
\begin{equation}\label{4.17}
A_{\la-1}<B_{\la-1}=N\quad \mbox{and}\quad B_\la=0.
\end{equation}

\medskip
Let us finally apply the change of coordinates $y_1:=x_1,\, y_2:=x_2-\rho(x_1).$  The non-zero roots $\tilde r$ of $\tilde f$ are then given by $\tilde r=r-\rho,$ with $r\in\dotol$ for some $l<\la,$  and they have the same multiplicities and leading exponents as  the corresponding roots $r.$ In view of \eqref{4.2}, the vertices of the Newton diagram $\N(\tilde f)$ are thus given by the points 
$(A_0,,B_0),\dots,(A_{\la-1},B_{\la-1}),$ i.e., the effect of the change of coordinates on the Newton diagram is the removal of the last vertex $(A_{\la},B_{\la}).$ Since, by \eqref{4.17}, the vertex  $(A_{\la-1},B_{\la-1})$  lies above the bisectrix, we see  that the principal face of  $\N(\tilde f)$ is the horizontal half-line emerging from the point $(A_{\la-1},B_{\la-1})$ along the line $t_2=N,$  a non-compact set.  According to Theorem \ref{s3.4},  the coordinates $y$ are thus adapted to $f,$  and the height $h(f)$ is given by 
\begin{equation}\label{4.17n}
h(f)=N.
\end{equation}

\qed

The proof suggests the following definitions. If
$$
r(x_1)=c_{l_1}^{\al_1}x_1^{a_{l_1}}+c_{l_1l_2}^{\al_1\al_2}x_1^{a_{l_1l_2}^{\al_1}}+\cdots+
c_{l_1\cdots l_p}^{\al_1\cdots \al_p}x_1^{a_{l_1\cdots
l_p}^{\al_1\cdots \al_{p-1}}}+\cdots,
$$ is any root of $f$ (more precisely, of $F$), then any leading part 
$$
\sum_{p=1}^K 
c_{l_1\cdots l_p}^{\al_1\cdots \al_p}x_1^{a_{l_1\cdots
l_p}^{\al_1\cdots \al_{p-1}}} ,
$$
with $1\le K\le\infty,$ which is real analytic, i.e, where all exponents appearing  in this series are positive integers, will be called an  {\it analytic root jet}. 

We have seen that the function $\psi$ constructed by Varchenko's algorithm in a unique way is indeed an analytic root jet, which we call the {\it principal root jet}.
\medskip

Our proof even reveals that the conditions (a) - (c) in Theorem \ref{s3.4} are necessary and sufficient for the adapteness of the given coordinate system for arbitrary analytic functions $f.$ In the statements of the following corollaries we shall always make the following 
\medskip 

\noi {\bf General Assumptions.} {\it The function $f\x$ is a real valued  real analytic function defined on a neighborhood of the origin in $\bR^2$ with
$f(0,0)=0,\,\nabla f(0,0)=0$. Choose $\ka_1,\ka_2\ge 0$ such that the principal face $\pi(f)$ of
the Newton polyhedron  of $f$  lies on the   line $\ka_1 t_1+\ka_2
t_2=1,$  and  assume that  $\ka_2\ge\ka_1$. }

\medskip
\begin{cor}\label{s4.3}
The given coordinates $\x$ are not adapted to $f$ if and only if  the conditions (a) - (c) in Theorem \ref{s3.4} are satisfied. In particular, the given coordinates are adapted to $f$ if and only if they are adapted to the principal part $f_p$ of $f.$ 

Moreover,  we always have $h(f)\le h(f_p).$
\end{cor}

\noi{\bf Proof.}  The necessity of these conditions for non-adaptedness has been proved in Theorem \ref{s3.4}. Assume conversely that (a) - (c) hold true. In that case, we have seen in the proof of Theorem \ref{s4.2} that there exists a change of coordinates which strictly increases the height, so that the original coordinates are not adapted. 

\medskip
Since the conditions (a) - (c) depend in fact only on $f_p,$ we see in particular that the given coordinates are adapted to $f$ if and only if they are adapted to the principal part $f_p$ of $f.$ 

\medskip
Finally, if the coordinates $\x$ are adapted to $f$ (hence also to $f_p$), then we clearly have $h(f)=h(f_p).$  Otherwise, the proof of Theorem \ref{s4.2} shows that the first change of coordinates $y_1:=x_1,\, y_2:= x_2-c_\la^\beta x_1^{a_\la}$ that we considered reduces the principal edge $\pi(f)= [(A_{\la-1},B_{\la-1}) ,(A_\la,B_\la)]$ to the shorter interval (possibly of length zero)  $[( A_{\la-1}, B_{\la-1}),(\tilde A_\la,\tilde B_\la)]$ on the same line, but lying above the bisectrix. Notice that the coordinates $(y_1,y_2)$ are already adapted to the principal part $f_p,$ and that, by  \eqref{4.14n}, \eqref{4.12} and Theorem \ref{s3.4}, we have $\tilde B_\la=N\bola=m(f_p)=h(f_p).$ 

But, the point $(\tilde A_\la,\tilde B_\la)$ will be contained in the Newton diagrams of $f$ associated to all subsequent systems of coordinates that we constructed by our algorithm (compare \eqref{4.14}, applied to the coordinates $y$), which means that the principal face of $f$ in our final, adapted coordinate system must lie in the half-space $t_2\le h(f_p), $ so that $h(f)\le h(f_p).$ 

\qed

\begin{cor}\label{s4.4}
a) We can always find a change of coordinates $x=\vp(y)$  at $0$ of the form $y_1:=x_1,\, y_2:=x_2-\psi(x_1),$  such that the coordinates  $(y_1,y_2)$ are adapted to  $\tilde f:=f\circ \vp$ and  the following hold true:

If the principal face $\pi(\tilde f)$ is compact and lies on the line $\tilde \ka_1 t_1+\tilde \ka_2 t_2=1, $ with $\tilde \ka_1\le \tilde \ka_2,$ then $\psi$ is a polynomial of degree strictly less than $\tilde \ka_2/\tilde \ka_1.$ 

This applies in particular, if the height $h(f)$ of $f$ is a non-integer  rational   number.

\medskip
b) There always exists  a change of coordinates $x=\vp(y)$ of the form $y_1:=x_1,\, y_2:=x_2-\eta(x_1)$ at the origin,  with a polynomial function   $\eta(y_1),$  such that in the new coordinates $y,$ we have 
 $h(f)=h(\tilde f)=h(\tilde f_p).$ Here, we have again put $\tilde f:=f\circ \vp.$
\end{cor}

\noi{\bf Proof.}  Indeed, the algorithm that we devised in the proof of Theorem \ref {s4.2} in order  to construct  an adapted coordinate system  shows that we can arrive at an adapted coordinate system with a polynomial function $\psi,$ unless we have to choose for $\psi$ one of the roots $r$ with infinitely many non-trivial terms in its  Puiseux series expansion. In the latter case, the principal face in the adapted coordinate system that we constructed is non-compact and the height is an integer, as we have seen.  Thus, if the principal face in the adapted coordinate system is compact, the algorithm must terminate after  a finite number of steps. And, in Step 1, the degree of the polynomial used in the change of coordinates is given by $a_\la,$ where by \eqref{4.10} $a_\la=\ka_2/\ka_1$ is just the inverse of the slope of the principal edge of the Newton polyhedron of $f.$  However, the proof shows that  the slope of the principal face strictly decreases by the change of coordinates in Step 1, and the same applies to all subsequent steps. If we apply this to the last change of coordinates before achieving adapted coordinates, we see that the function $\psi$ is a polynomial of degree $m<\tilde \ka_2/\tilde \ka_1.$ 
 This proves a).

Moreover, in the  case where our algorithm does not terminate, after finitely steps in our algorithm (which all consist of polynomial changes of coordinates), we may assume that the  polynomial $\tilde f_p$ corresponding to the  principal face of the Newton diagram has a unique root of multiplicity $N=h(f)$ (compare \eqref{4.15}, \eqref{4.16} and \eqref{4.17n}).  
If we choose the coordinates $y$ which we obtain at this stage, we then have $h(\tilde f_p)=m(\tilde f_p)=N=h(f)=h(\tilde f),$  so that also b) is proven. 

\qed

\bigskip

\setcounter{equation}{0}
\section{The smooth case}\label{smooth}

We shall finally extend Theorem \ref{4.2} to the smooth setting.

\begin{thm} \label{s5.1}
Let $f$ be a real valued smooth  function of finite type defined on a neighborhood of the origin in $\bR^2$ with
$f(0,0)=0,\,\nabla f(0,0)=0$. Choose $\ka_1,\ka_2\ge 0$ such that the principal face $\pi(f)$ of
the Newton polyhedron  of $f$  lies on the   line $\ka_1 t_1+\ka_2
t_2=1.$  Without loss of generality,  we may assume that  $\ka_2\ge\ka_1$. Then there exists a  smooth function  $\psi(x_1)$ of $x_1$ near the origin with $\psi(0)=0$ such that an adapted coordinate
system $(y_1,y_2)$ for $f$ near $0$ is given by  $y_1:=x_1,\, y_2:=x_2-\psi(x_1).$ 
\end{thm}

\medskip

\noi{\bf Proof.}  We proceed in a very similar way as in the analytic setting, adopting also the same notation.  If the coordinates are adapted to $f,$  then we may choose $\psi:=0$ and are finished. 
\medskip

Otherwise,  again by Theorem \ref{s3.4},   the principal face $\pi(f)$ is a  compact edge. Moreover, we have $ \ka_2/\ka_1=:m_1\in \bN,$ and $m(f_p)=m(f_\ka)>d(f).$ Let us then choose a real root $x\mapsto b_1x_1^{m_1}$ of the principal part $f_p$ of $f$ of maximal multiplicity $N_0:=m(f_p),$ i.e., the principal root. Then $b_1\ne 0, $ again by Theorem \ref{s3.4},  

\medskip
\noi{\bf Step 1.} We apply the  real  change of variables $x=\vp(y)$ given by $y_1:= x_1, \, y_2:= x_2-b_1x_1^{m_1},$  and put  $\tilde f:=f\circ \vp.$  Let us again endow all quantities associated to $\tilde f$ with a superscript $\widetilde{}$.  Now,  if the coordinates $y$ are adapted to $f,$ we choose $\psi(x_1):= b_1x_1^{m_1}$ and are finished. 
\medskip

Otherwise, the principal face $\pi(\tilde f)$ is a  compact edge, and  we have $ \tilde \ka_2/\tilde \ka_1=:m_2\in \bN$ and $N_1:=m(\tilde f_p)>d(\tilde f).$ Recall that the principal part $\tilde f_p$ of $\tilde f$ is $\tilde \ka$-homogeneous of degree one. 
We claim that 
\begin{equation}\label{5.1}
m_2>m_1 \quad \mbox{and} \quad N_1\le N_0.
\end{equation}

Indeed, recall that  the effect of our change of coordinates $\vp$ on the Newton polyhedron is such that it preserves all lines $\ka_1t_1+\ka_2 t_2=c. $ Let us therefore 
choose $m\in\bN $ so big that we have $j+k\le m$ for every point $(j,k)$ lying on any  such line  passing through any point in the Newton diagram $\N_d(\tilde f)$ of $\tilde f,$ and denote by $F$ the Taylor polynomial of order $m$ of $f.$ Then it is clear that $f$ and $F$ have the same principal faces and parts, and the same applies to  $\tilde f$ and $\tilde F:=F\circ \vp.$ I.e., we have $f_p=F_p$ and $\tilde f_p=\tilde F_p.$ We can therefore apply our results for the analytic case to the polynomial function $F$ and obtain \eqref{5.1}.  
\smallskip

This argument also shows that the change of coordinates increases the distance, i.e., $d(\tilde f)> d(f).$ 

\medskip
\noi{\bf Subsequent steps.} 
Now, either the new coordinates $y$ are adapted, in which case we are finished. Or we can apply the same procedure to $\tilde f,$   etc.. In this way, we obtain a sequence of functions $f_{(k)}=f\circ \vp_{(k)}, $ with $f_{(0)}:=f$ and $f_{(k+1)}:=\widetilde{f_{(k)}},$  which can be obtained from the original coordinates $x$ by means of a change of coordinates $x=\vp_{(k)}(y)$ of the form 
$$y_1:=x_1,  \, y_2:= x_2-\sum_{l=1}^{k}b_lx_1^{m_l},
$$
with positive integers $m_1<m_2<\cdots <m_k<m_{k+1}<\cdots$  and real coefficients $b_l\ne 0.$ Moreover, if $N_k:=m((f_{(k)})_p)$ denotes the maximal order of vanishing of the principal part of $f_{(k)}$ along the unit circle $S^1,$ then we have 
$$N_0\ge N_1\ge \cdots\ge N_{k}\ge N_{k+1}\ge \cdots.$$

Either  this procedure will stop after finitely many steps, or it will continue infinitely. If it stops, say, at the $k$-th step, it is clear that we will have arrived at an adapted coordinate system $x=\vp_{(k)} (y), $ with a polynomial function $\psi(x_1)=\sum_{l=1}^{k}b_lx_1^{m_l}.$

\medskip

\noi{\bf Final step.} Assume that the procedure does not terminate. Since the maximal  multiplicities $N_k$ of the roots of  the principal part of $f_{(k)}$ form a decreasing sequence, we find again some $k_0, N\in \bN$ such that $N_k=N$ for every $k\ge k_0.$ By comparing with the effect of the change of coordinates in each step of order $k\ge k_0$ with the effect on the Taylor polynomial of sufficiently high degree, we see from the corresponding result \eqref{4.17a} in the analytic case that the principal part of $f_{(k)}$ is of the form
$$
(f_{(k)})_p(x)=c_kx_1^{\nu_1} (x_2-b_{k+1} x_1^{m_{k+1}})^N,
$$
where $\nu_1<N.$ 
It is mixed homogeneous of degree one with respect to the weight $\ka^{(k)},$ given by 
$$\ka^{(k)}_1:= \frac 1{\nu_1+Nm_{k+1}},\quad  \ka^{(k)}_2:= \frac {m_{k+1}}{\nu_1+Nm_{k+1}},
$$
so that, for $k\ge k_0,$ 
\begin{equation}\label{5.2}
f_{(k)}(x)=c_kx_1^{\nu_1} (x_2-b_{k+1} x_1^{m_{k+1}})^N \ +\ \mbox{terms of higher $\ka^{(k)}$- degree}.
\end{equation}
\medskip

Now, according to the classical Borel lemma, we can find a smooth function $\rho(x_1)$ near the origin whose Taylor series is the formal series $\sum_{l=1}^{\infty}b_lx_1^{m_l}.$ Consider the smooth change of coordinates  $x:=\vp(y)$ given by  $y_1:=x_1,\, y_2:=x_2-\rho(x_1),$ and put $\tilde f:=f\circ \vp.$ We claim that the coordinates $y$ are adapted to $\tilde f.$ 

\medskip
Indeed,  we have 
$$
\tilde f(y)= f_{(k)}\circ (\vp_{(k)}^{-1}\circ \vp)(y)=f_{(k)}\Big(y_1,y_2+(\rho(y_1)-\sum_{l=1}^{k}b_ly_1^{m_l})\Big),
$$
where $\rho(y_1)-\sum_{l=1}^{k}b_ly_1^{m_l}$ has the Taylor series $\sum_{l=k+1}^{\infty}b_ly_1^{m_l}.$ In view of  \eqref{5.2}, this show that 
\begin{equation}\label{5.3}
\tilde f(y)=c_ky_1^{\nu_1} y_2^N \ + R(y),
\end{equation}
where $R$ is a smooth function consisting of terms of $\ka^{(k)}$-degree strictly bigger than one, for every $k\ge k_0,$ i.e., 
\begin{equation}\label{5.4}
 \frac 1{\nu_1+Nm_{k+1}} j_1+   \frac {m_{k+1}}{\nu_1+Nm_{k+1}} j_2>1 \quad \mbox{for every}\ (j_1,j_2)\in \T(R).
\end{equation}
Since $m_k\to\infty$ as $k\to \infty,$ this implies that $j_2\ge N.$ Moreover, if $j_2=N, $ then the left-hand side of \eqref{5.4} is bounded by $1,$ if $j_1\le \nu_1,$ so that we must have $j_1>\nu_1.$ In combination, \eqref{5.3} and \eqref{5.4} show that the Newton polyhedron $\N(\tilde f)$ of $\tilde f$ contains the point $(\nu_1,N),$ but nor further point on the left to this point on the line $t_2=N,$ and that all other points of $\N(\tilde f)$ are contained in the open half-plane above this line. Since $\nu_1<N, $ this shows that the principal face of $\N(\tilde f)$ is the unbounded horizontal half-line given by $t_1\ge \nu_1,\ t_2=N,$ hence the coordinates are adapted to $\tilde f.$

\qed

\medskip

The following corollary is immediate from the proof  of Theorem \ref{s5.1}  and the proofs of  the Corollaries \ref{s4.3} and \ref{s4.4}, which carry over to the smooth setting. 

\begin{cor}\label{s5.2}
Let $f$ be a real valued  smooth  function of finite type defined on a neighborhood of the origin in $\bR^2$ with
$f(0,0)=0,\,\nabla f(0,0)=0$. Then the height  $h(f)$ is a rational number. 

Moreover, if we choose $\ka_1,\ka_2\ge 0$ such that the principal face $\pi(f)$ of
the Newton polyhedron  of $f$  lies on the   line $\ka_1 t_1+\ka_2
t_2=1,$  and assume without loss of generality that  $\ka_2\ge\ka_1,$ then the  
Corollaries \ref{s4.3} and \ref{s4.4}  remain true in this smooth setting. 
\end{cor}

\bigskip

\begin{center} {ACKNOWLEDGEMENT}
\end{center}
We wish to thank Michael Kempe for his   support in creating the graphics included in this article.

\end{document}